\documentclass[matprg,final]{svjour}

\usepackage{amsmath}
\usepackage{graphicx}
\usepackage{amsfonts}
\usepackage{amssymb}
\usepackage{url}

\newcommand{\re}{\mathbb{R}}

\newcommand{\N}{\mathbb{N}}

\newcommand{\lmd}{\lambda}
\newcommand{\half}{\frac{1}{2}}

\newcommand{\nn}{\nonumber}
\newcommand{\eps}{\epsilon}
\newcommand{\veps}{\varepsilon}
\newcommand{\vareps}{\varepsilon}
\newcommand{\dt}{\delta}
\newcommand{\Dt}{\Delta}
\newcommand{\af}{\alpha}

\newcommand{\sig}{\sigma}

\newcommand{\supp}{\mbox{supp}}

\newcommand{\eproof}{{\ $\square$ }}
\newcommand{\reff}[1]{(\ref{#1})}
\newcommand{\pt}{\partial}
\newcommand{\prm}{\prime}

\newcommand{\mc}[1]{\mathcal{#1}}

\def\cL{ {\mathcal L} }

\def\cM{{\mathcal M}}

\def\bdS{\partial S}
\def\tG{\tilde G}
\def\tS{\tilde S}
\def\hS{{\hat S}}
\def\tG{ { \tilde G}}

 \def\bdS{\partial S}
 \def\bdiS{ \partial S_i \cap S}

\newcommand{\bdes}{\begin{description}}
\newcommand{\edes}{\end{description}}

\newcommand{\bal}{\begin{align}}
\newcommand{\eal}{\end{align}}

\newcommand{\bnum}{\begin{enumerate}}
\newcommand{\enum}{\end{enumerate}}

\newcommand{\bit}{\begin{itemize}}
\newcommand{\eit}{\end{itemize}}

\newcommand{\bea}{\begin{eqnarray}}
\newcommand{\eea}{\end{eqnarray}}

\newcommand{\be}{\begin{equation}}
\newcommand{\ee}{\end{equation}}

\newcommand{\baray}{\begin{array}}
\newcommand{\earay}{\end{array}}

\newcommand{\bsry}{\begin{subarray}}
\newcommand{\esry}{\end{subarray}}

\newcommand{\bca}{\begin{cases}}
\newcommand{\eca}{\end{cases}}

\newcommand{\bcen}{\begin{center}}
\newcommand{\ecen}{\end{center}}

\newcommand{\bbm}{\begin{bmatrix}}
\newcommand{\ebm}{\end{bmatrix}}

\newcommand{\bmx}{\begin{matrix}}
\newcommand{\emx}{\end{matrix}}

\newcommand{\bpm}{\begin{pmatrix}}
\newcommand{\epm}{\end{pmatrix}}

\newcommand{\btab}{\begin{tabular}}
\newcommand{\etab}{\end{tabular}}

\newtheorem{prop}[theorem]{Proposition}
\newtheorem{lem}[theorem]{Lemma}
\newtheorem{cor}[theorem]{Corollary}

\newtheorem{exm}[theorem]{Example}

\renewcommand{\subsection}[1]{
     \stepcounter{subsection}
     \settowidth{\hangindent}{\bf\thesubsection.~}
     \bigskip\noindent
     \noindent
     {\bf\hbox{\thesubsection.~}#1}\par
     \nobreak
     \medskip
}

\usepackage[top=1.5in,bottom=1.2in,left=1.1in,right=1.1in]{geometry}

\begin{document}

\title{Semidefinite Representation of Convex Sets}

\author{J. William Helton \and Jiawang Nie}

\institute{
J. William Helton \at
Department of Mathematics, University of California at San Diego,
9500 Gilman Drive, La Jolla, CA 92093.
\email{helton@math.ucsd.edu} \and
Jiawang Nie \at
Department of Mathematics, University of California at San Diego,
9500 Gilman Drive, La Jolla, CA 92093.
\email{njw@math.ucsd.edu}
}

\maketitle



\begin{abstract}
Let $S =\{x\in \re^n:\, g_1(x)\geq 0, \cdots, g_m(x)\geq 0\}$ be a
semialgebraic set defined by multivariate polynomials $g_i(x)$.
Assume $S$ is convex, compact and has nonempty interior. Let $S_i
=\{x\in \re^n:\, g_i(x)\geq 0\}$,
and $\bdS$ (resp. $\bdS_i$) be the boundary of $S$ (resp. $S_i$).
This paper, as does the subject of
semidefinite programming (SDP), concerns Linear Matrix
Inequalities (LMIs). The set $S$ is said to have an LMI
representation if it equals the set of solutions to some LMI and
it is known that some convex $S$ may not be LMI representable
\cite{HV}. A question arising from \cite{NN94}, see \cite{HV,N06},
is: given a subset $S$ of $\re^n$, does there exist an LMI representable
set $\hS$ in some higher dimensional space $ \re^{n+N}$ whose
projection down onto $\re^n$ equals $S$. Such  $S$ is called
semidefinite representable or SDP representable. This paper
addresses the SDP representability problem.

\smallskip

The following are the main contributions of this paper: {\bf (i)}
Assume $g_i(x)$ are all concave on $S$. If the positive definite
Lagrange Hessian (PDLH) condition holds, i.e., the Hessian of the
Lagrange function for optimization problem of minimizing any
nonzero linear function $\ell^Tx$ on $S$ is positive definite at
the minimizer, then $S$ is SDP representable. {\bf (ii)} If each
$g_i(x)$ is either sos-concave ($-\nabla^2g_i(x)=W(x)^TW(x)$
for some possibly nonsquare matrix polynomial $W(x)$)
or strictly quasi-concave on $S$,
then $S$ is SDP representable.
{\bf (iii)} If each $S_i$ is either sos-convex
or poscurv-convex
($S_i$ is compact convex, whose boundary has  positive curvature
 and is nonsingular, i.e. $\nabla g_i(x) \not = 0$ on $\bdS_i \cap S$),
then $S$ is SDP representable.
This also holds for $S_i$ for which
$\bdS_i \cap S$ extends smoothly to the boundary of
a poscurv-convex set containing $S$. {\bf (iv)} We give the
complexity of Schm\"{u}dgen and Putinar's matrix
Positivstellensatz, which are critical to the proofs of (i)-(iii).
\end{abstract}

\bigskip


\section{Introduction}
One of the main advances in optimization which has had a profound
effect on control theory and nonconvex optimization as well as
many other disciplines is semidefinite programming (SDP)
\cite{N06,SDPbook}. This gives effective numerical algorithms
for solving problems presented in terms of Linear Matrix
Inequalities (LMIs). Arising from this is the very basic issue of
which problems can be presented with LMIs and this paper addresses
one of  the most classical aspects of this problem.

We say a set $S$ have an {\it LMI representation} or is {\it
LMI representable} if
\be
S=\{x\in \re^n: A_0+\sum_{i=1}^n A_i x_i \succeq 0\}
\ee
for some symmetric matrices $A_i$. Here the notation $X
\succeq 0 \,(\succ 0)$ means the matrix $X$ is positive
semidefinite (definite). If $S$ has an interior point, $A_0$ can
be assumed to be positive definite without loss of generality.
Obvious necessary conditions for $S$ to be LMI representable are
that $S$ must be convex and $S$ must also be a basic closed semialgebraic
set
\[
S =\{x\in \re^n:\, g_1(x)\geq 0, \cdots, g_m(x)\geq 0\}
\]
where $g_i(x)$ are multivariate polynomials.
We shall always assume $S$ has an interior point. For example, any
convex quadratic constraint $ \{ x \in \re^n:\,
a+b^Tx-x^TC^TCx\geq 0\}$ can be represented by the LMI
\[
\left\{ x\in \re^n:\, \bbm I_n & Cx \\
(Cx)^T & a+b^Tx \ebm \succeq 0 \right\}
\]
where $I_n$ is the $n\times n$ identity matrix. Here $B^T$ denotes
the transpose of matrix $B$. A basic question (asked in \cite{PS})
is: which convex sets can be represented by LMIs? It turns out
that some convex sets are not LMI representable. Helton and
Vinnikov \cite{HV} proved that a strong condition called {\it
rigid convexity} is  necessary  for a set to have an LMI
representation (as well as sufficient in case of dimension two).
For instance, the convex set
\be
\label{eq:TV}
T = \left\{x\in \re^2:\, 1-(x_1^4+x_2^4) \geq 0 \right\}
\ee
does not admit an LMI representation \cite{HV}, since it is not
rigidly convex.

However, the set $T$ is the projection onto $x$-space of the set
$$
\hS:= \left\{(x,w)\in \re^2\times \re^2:\, \bbm 1+w_1 & w_2 \\ w_2
& 1-w_1 \ebm \succeq 0, \bbm 1 & x_1 \\ x_1 & w_1 \ebm \succeq 0,
\bbm 1 & x_2 \\ x_2 & w_2 \ebm \succeq 0 \right\}
$$
in $\re^4$ which is represented by  an LMI. This motivates
\medskip

\noindent {\bf Question:} \ {\it Which convex sets $S$ are the
projection of a set $\hS$ having an LMI representation; in other
words, do there exist symmetric
  matrices $F_i,G_j$ such that
$S$ equals $\{ x : \ (x,y) \in \hS \}$ where
\be \label{eq:liftLMI}
\hS =\left\{ (x, y) \in \re^{(n+N)}: F_0 +
\sum_{i=1}^n F_i x_i + \sum_{j=1}^N G_j y_j \succeq 0\, \right\} .
\ee }

\noindent Such sets $S$ are called {\it semidefinite
representable} or {\it SDP representable}.
Ben-Tal and Nemirovskii (\cite{BN}),
Nesterov and Nemirovskii
(\cite{NN94}), and Nemirovskii (\cite{N06})
gave collections of examples of SDP representable sets. Thereby
leading  to the question  which sets are SDP representable? In \S
4.3.1 of his excellent 2006 survey \cite{N06} Nemirovsky said ``
this question seems to be completely open". Obviously, to be SDP
representable, $S$ must be convex and semialgebraic. What are the
sufficient conditions that guarantee $S$ is SDP representable?
This paper
addresses this kind of question. Sometimes we refer to a
semidefinite representation as a {\it lifted LMI representation}
of the convex set $S$ and to the LMI in \reff{eq:liftLMI} as the
{\it  lifted LMI for $S$}.

\medskip

A construction of the SDP representation for convex sets was
proposed  by Lasserre \cite{Las06} and also in
the dimension two case by
Parrilo, for example in \cite{Par06},
and could be viewed using
the following idea. Let $\cM $ denote the
space of Borel measures on $S$ and let $\hS$ denote the convex
subset of all nonnegative mass one measures. The Krein Millman
Theorem \cite{C} says that $\hS$ projects down onto $S$ via
\[
P(\mu) := \int_S x d\mu(x) \qquad \quad \mu \in \hS.
\]
Unfortunately $\hS$ is infinite dimensional, so unsuitable as an
SDP representation.
The Lasserre and Parrilo proposal, which will be sketched  later,
is to cut down $\hS$ by looking at it as the set of all positive
mass one linear functionals on the polynomials of some fixed
degree $N$. Moment and sum of squares (SOS) techniques show that
this gives an LMI, denoted by $\cL_N$, for each degree $N$, and
that the projection onto $x-$ space of the set $\hS_N: = \{(x, y):
\mathcal{L}_N(x,y) \geq 0 \}$ contains $S$ for all $N$. The open
question remaining is whether there exists an integer $N$ large
enough to produce the equality.

The validity  of this general type of  construction has been
supported by very nice recent findings on the SDP  representation
of convex sets.
Parrilo \cite{Par06} proved this gives a lifted
LMI representation in the two dimensional case when the boundary
of $S$ is a single rational planar curve of genus zero. Lasserre
\cite{Las06} proved this construction can give arbitrarily accurate
approximations when $N$ goes to infinity.

This article gives sufficient conditions
(presented as the
hypotheses of Theorems \ref{thm:conPDLH} through \ref{thm:posCurv})
on a convex set $S$ guaranteeing that it is SDP representable.
The first condition we present (Theorem \ref{thm:conPDLH})
is the PDLH condition, and we prove validatity of
Lasserre-Parrilo moment type constructions when PDLH holds.
After that  come three
theorems, each having
weaker hypotheses than
the preceding one
and each concluding that $S$ is SDP representable.
The last of them (Theorem \ref{thm:posCurv})
is a bit weaker than saying
that a strictly convex basic semialgebraic set $S$
has an SDP representation,
provided its boundary is nonsingular
(the gradients of defining polynomials for $S$
do not vanish).

More specifically,
Theorems \ref{thm:setSimple} and \ref{thm:posCurv}
are based on geometric properties of the boundary
of a convex set. Any  convex set has boundary,
which if a smooth manifold,
has nonnegative curvature
and conversely,
having positive curvature everywhere is slightly stronger
than strict convexity.
Strict convexity and positively curved boundary
are not the same as is illustrated by
the set in example \eqref{eq:TV}
which while strictly convex has zero curvature at
$(-1,0)$, $(1,0)$, $(0,-1)$ and $(0,1)$,
although it has positive curvature everywhere else.
A good illustration of the gap between our necessary and sufficient
conditions for SDP representability is Theorem \ref{thm:setSimple}
when specialized to a convex set $S$
defined by a single polynomial $g(x)$.
It implies that if
the boundary $\bdS$ is nonsingular
($\nabla g(x) \not = 0$ on $\bdS$)
and if $\bdS$ has positive curvature at all points,
then $S$ has an SDP representation.
Thus for a single defining function $g(x)$
the necessary vs. sufficient gap lies only in
the gradient being zero or the curvature being zero
somewhere on $\bdS$.
Theorems \ref{thm:setSimple} and
\ref{thm:posCurv} also give generalizations of this
for more than one defining function.


A subsequent paper, based on this one,  \cite[Section~3]{HN2},
extends these results and shows that if each component
of the boundary of convex and bounded
$S$ is positively curved and nonsingular,
then $S$ is SDP representable.

We should emphasize that while our description stresses
a clean characterization of
existence of LMI lifts, we shall soon introduce a class of sets
we call SOS- convex and constructions for them which might
be practical on modest size problems.
Now we turn to a formal presentation of results.


Let $S =\{ x\in \re^n:\, g_1(x)\geq 0,\, \cdots,\, g_m(x)\geq 0
\}$ be a basic closed semialgebraic set; here the $g_i(x)$ are in
the ring $\re[x]$ of multivariate polynomials with real
coefficients and are called the {\it defining polynomials} for
$S$. Assume $S$ is convex, compact and has nonempty interior.
Let
$S_i =\{ x\in \re^n:\, g_i(x)\geq 0 \}$ and $Z(g_i) =\{ x\in
\re^n:\, g_i(x) = 0 \}$ be the zero set of $g_i$.
Denote by $\bdS$ and $\bdS_i$ the
boundaries of $S$ and $S_i$ respectively.
Note that $\bdS_i$ might be contained in $Z(g_i)$ properly.

\medskip

First, consider the case that all the defining polynomials
$g_i(x)$ are concave on $S$. The {\it positive definite Lagrange
Hessian (PDLH)} condition requires that for any nonzero vector
$\ell\in\re^n$, the Hessian of the Lagrange function
corresponding to the optimization problem of minimizing $\ell^Tx$
over $S$ is positive definite at each minimizer, i.e.,
$-\sum_{i=1}^m \lmd_i \nabla^2 g_i(u)$ is positive definite for
every minimizer $u$ and the corresponding Lagrange
multipliers $\lmd_i\geq 0$.
Obviously, if every $g_i(x)$ has negative definite Hessian
on the boundary, then the PDLH condition holds.

\begin{theorem} \label{thm:conPDLH}
Suppose $S =\{x\in \re^n:\, g_1(x)\geq 0, \cdots, g_m(x)\geq 0\}$
is compact convex and has nonempty interior.
Assume $g_i(x)$ are concave
on $S$. If the PDLH condition holds, then $S$ is SDP
representable.
\end{theorem}

\noindent
{\it Remark:} In Theorem~\ref{thm:conPDLH}, where the
$g_i(x)$ are concave on $S$, the matrix $-\sum_{i=1}^m \lmd_i
\nabla^2 g_i(u)$ must be positive semidefinite. The PDLH condition
requires it is positive definite, i.e., its determinant is
nonzero, which defines a Zariski open set. So the PDLH condition is a
generic condition subject to the property that $g_i(x)$ are
concave on $S$.

The SDP representation of $S$ in Theorem~\ref{thm:conPDLH}
can be constructed explicitly,
which will be shown in Section~\ref{sec:LMIlift}.
The
lifted LMI~\reff{LmiN} or \reff{putLMI}
(archimedean condition is then required)
represents $S$ exactly under the PDLH condition.
However, when the PDLH condition fails,
the constructed LMIs in Section~\ref{sec:LMIlift}
might not represent $S$ correctly.
This leads to our next result.

\medskip
Second,
consider the case that all the defining polynomials
$g_i(x)$ are {\it quasi-concave} on $S$.
This means the super
level set $S_i(\af)=\{x\in S:\, g_i(x)\geq \af\}$ is convex
for every $\af \in g_i(S)$.
So the level set
$Z(g_i-\af)=\{x\in \re^n:\, g_i(x)= \af\}$ when smooth has {\it nonnegative
curvature} in $S$, i.e., for all $x\in S$,
\begin{align*}
-v^T\nabla^2 g_i(x)v \geq 0 ,
 \ \ \ \forall\, v \in \nabla g_i(x)^\perp =
\{ v\in \re^n:\, v^T \nabla g_i(x) =0\}.
\end{align*}
We say $g_i(x)$ is {\it strictly quasi-concave} on $S$ if every
$Z(g_i-\af)$ has {\it positive curvature} in $S$, i.e.,
\begin{align*}
-v^T\nabla^2 g_i(x)v > 0 , \  \  \  \forall\, 0\ne v \in \nabla
g_i(x)^\perp, \, \forall\, x\in S.
\end{align*}
By
Exercise~3.44(a) in Boyd and Vandenberghe \cite{BV}, the above
is equivalent to the {\it modified Hessian}
\begin{align*}
-\nabla^2 g_i(x) + M \nabla g_i(x) \nabla g_i(x)^T \succ 0
\end{align*}
for some constant $M>0$. It will be shown (Lemma~\ref{modHes})
that the constant $M$ can be chosen uniformly for all $x\in S$ if
$g_i(x)$ is strictly quasi-concave on $S$.

A polynomial $p(x)$ is said to be a {\it sum of squares (SOS)} if
$p(x) = w(x)^Tw(x)$ for some column vector polynomial $w(x)$. The
necessary condition for $p(x)$ to be SOS is that it is nonnegative
on the whole space $\re^n$, but the converse might not be true. We
refer to \cite{Rez00} for a survey on SOS polynomials. A symmetric
matrix polynomial $P(x) \in \re[x]^{n\times n}$ is SOS if there
exists a possibly nonsquare matrix polynomial $W(x)$ with $n$ columns
such that $P(x)=W(x)^TW(x)$.
The defining polynomial $g_i(x)$ is called {\it
sos-concave} if the negative Hessian $-\nabla^2 g_i(x)=
-\left(\frac{\pt^2 g_i}{\pt x_k \pt x_\ell}\right)$ is SOS.
Similarly, $g_i(x)$ is called {\it sos-convex} if the Hessian
$\nabla^2 g_i(x)$ is SOS.

\begin{theorem} \label{thm:mainRep}
Suppose $S =\{x\in \re^n:\, g_1(x)\geq 0, \cdots, g_m(x)\geq 0\}$
is compact convex and has nonempty interior. If each $g_i(x)$ is either
sos-concave or strictly quasi-concave on $S$, then $S$ is SDP
representable.
\end{theorem}

\noindent
{\it Remark:} The case where a piece of the boundary is
a linear subspace  is included in Theorem~\ref{thm:mainRep},
 since if some $g_i(x)$ is a linear polynomial,
 then its Hessian is identically zero and hence $g_i(x)$ is sos-concave.

It is possible that the defining functions of a convex set
can be neither concave nor quasi-concave, because
the defining functions of a convex set
can behave badly in the interior.
However, they have nice properties near the boundary
which are helpful for us to establish
the semidefinite representability.
This leads to the following results.

\medskip
Third, consider the case that $S$ is convex but the defining
polynomials $g_i(x)$ are not quasi-concave on $S$. This is because
the super level sets of $g_i(x)$ might not be all convex. We call
$S_i$ {\it poscurv-convex} if $S_i$
is compact convex, its boundary $\bdS_i$ equals $Z(g_i)$,
and $\bdS_i$ is nonsingular (the gradient does not vanish)
and has positive curvature at each point on it, which means that
\be
\label{def:posCurvature} -v^T\nabla^2 g_i(x)v > 0 , \  \ \forall\,
0\ne v \in \nabla g_i(x)^\perp , \,  \forall \, x\in \bdS_i.
\ee
Note that the definition of poscurv-convex sets
also applies to  $g_i(x)$ which are {\it smooth functions}
(not necessarily polynomials).
$S_i$ is called {\it sos-convex}
if $g_i$ is a polynomial and $g_i$ is sos-concave.

\begin{theorem} \label{thm:setSimple}
Suppose $S =\{x\in \re^n:\, g_1(x)\geq 0, \cdots, g_m(x)\geq 0\}$
is compact convex and has nonempty interior.
If each $S_i$ is
either sos-convex or poscurv-convex, then $S$ is SDP
representable.
\end{theorem}

We turn now to more general sets $S_i$.
We say $S_i$ is {\it extendable poscurv-convex with respect to} $S$
if $g_i(x)>0$ whenever $x\in S$ lies in the interior of $S_i$ and
there exists a
poscurv-convex set $T_i=\{x: f_i(x) \geq 0\} \supseteq S$
such that $\pt T_i\cap S = \bdS_i \cap S$.
Here $f_i(x)$ is a smooth function (not necessarily a polynomial)
such that $T_i$ is compact convex, $\pt T_i = Z(f_i)$ and
$\pt T_i$ is nonsingular.
In other words, $\bdS_i\cap S$ can be extended to
become a part of the boundary of a poscurv-convex set
defined by a smooth function.

\begin{theorem} \label{thm:posCurv}
Suppose $S =\{x\in \re^n:\, g_1(x)\geq 0, \cdots, g_m(x)\geq 0\}$
is compact convex and has nonempty interior.
If each $S_i$ is
either sos-convex or extendable poscurv-convex with respect to
$S$, then $S$ is SDP representable.
\end{theorem}

Theorem~\ref{thm:setSimple} is a special case of
Theorem~\ref{thm:posCurv},
since every poscurv-convex set is of course extendable poscurv-convex.
However, it turns out, see
the follow-up paper \cite{HN2} to this one,
that extendable poscurv-convexity of $S_i$ with respect to $S$
does not require much more
than  the boundary $\bdS_i \cap \bdS$ has positive curvature.
In \cite{HN2} this is combined with Theorem \ref{thm:posCurv}
to obtain a stronger result:
if for every $i$ either $S_i$ is sos-convex or $\bdS_i \cap \bdS$ is
positively curved and nonsingular, then $S$ is SDP representable.

The proofs for the above theorems are based  on a variety of
techniques, and they produce new results which might be of
interest independent of SDP representation.
First, the proofs introduce a natural
technique of writing a polynomial as a sum of squares by twice
integrating its Hessian, which is very suited to handling
sos-concavity.
Second, we give degree bounds for
polynomials appearing in Schm\"{u}dgen's and Putinar's matrix
Positivstellensatz, see the Appendix \S \ref{sec:complexityPosSS}.
These two techniques allow us to obtain
bounds on the degrees of polynomials which appear in SOS
representations.
Third, it is possible that the set
$S_i=\{x\in\re^n:\, g_i(x)\geq 0\}$ is strictly convex but that
the polynomial $g_i(x)$ is neither concave nor quasi-concave.
In \S \ref{sec:ConDef} we show under modest hypotheses that there
is a new set of (local) defining polynomials $p_i$ for the $S_i$
which are  strictly concave. This allows us to prove
Theorem~\ref{thm:posCurv} by using the new defining polynomials
$p_i$ together with the original polynomials $g_i$.


\smallskip
Now we say a few words about constructions.
When all the defining polynomials $g_i(x)$ are sos-concave,
an explicit SDP representation for $S$ is given by \reff{lassLMI}
in Section~\ref{sec:globconcav}. The forthcoming article
\cite{HN3} illustrates this in several examples,
and also shows how to incorporate sparsity.
If all $g_i(x)$ are  concave on $S$
and every $g_i(x)$ is either sos-concave or strictly concave
(having negative definite Hessian)
on the boundary $\bdS$ where it vanishes,
an explicit SDP representation for $S$ is given
by \reff{LmiN} or \reff{putLMI}
in Section~\ref{sec:LMIlift}
when the relaxation order $N$ is big enough.
For the time being, we do not have an estimate of how
large $N$ is sufficient.
In Theorem~\ref{thm:conPDLH},
the SDP representation can be constructed in the same way
as in Section~\ref{sec:LMIlift}.
In Theorems~\ref{thm:mainRep}, \ref{thm:setSimple} and \ref{thm:posCurv},
we have only shown the existence of SDP representations for $S$.
We would expect  it would
be very difficult to use the proof there  constructively.




The following notations will be used. For $x\in\re^n$,
$\|x\|=\sqrt{\sum_{i=1}^n x_i^2}$. $\N$ denotes the set of
nonnegative integers. For $\af \in \N^n$, $|\af| :=
\af_1+\cdots+\af_n$, $x^\af := x_1^{\af_1}\cdots x_n^{\af_n}$.
A $d$-form is a homogenous polynomial of degree $d$.
Given a set $K\subset \re^n$, $C^k(K)$ denotes the set of
$k$-times continuously differentiable functions in an open set
containing $K$, and $C^\infty(K)$ denotes the set of infinitely
times differentiable (smooth) functions in an open set containing
$K$. Given $\af\in\N^n$ and $f(x)\in C^k(K)$, $D^\af f(x) :=
\frac{\pt^{|\af|} f(x)} {\pt x_1^{\af_1}\cdots \pt x_n^{\af_n}}$.
For a symmetric matrix $A$, $\lmd_{\min}(A)$ denotes the smallest
eigenvalue of $A$, and $\| A\|_2$ denotes the standard $2$-norm of $A$.

The paper is organized as follows. Section~\ref{sec:LMIlift} gives
the constructions of lifted LMIs when $g_i$ are concave on $S$,
and states some theorems about
the sharpness of these lifted LMIs, whose
proofs will be given in Section~\ref{sec:proof}. Then
Section~\ref{sec:globconcav} turns to Lasserre and Parrilo
moment type of
constructions of lifted LMIs in \cite{Las06,Par06}, and we give a
sufficient condition that guarantees these constructed LMIs are
the SDP representations of $S$.
Section~\ref{sec:ConDef} discusses how to find concave defining
functions for poscurv-convex sets used to prove
Theorem~\ref{thm:posCurv}.
Section~\ref{sec:proof} gives proofs of the
theorems in Section~\ref{sec:LMIlift} and
in the Introduction.
Section \ref{sec:complexityPosSS} is an appendix
bounding the degrees of polynomials arising
in Schm\"{u}dgen's and Putinar's matrix Positivstellensatz.
Section~\ref{sec:cncl} summarizes
conclusions of the paper.

\section{The SDP representations when $g_i(x)$ are concave on $S$ }
\label{sec:LMIlift} \setcounter{equation}{0}

In this section we assume $g_i(x)$ are concave on $S$. Two kinds
of SDP representations will be constructed.
The construction of these SDP representations which
we use can be found in Lasserre \cite{Las06}.
We review the construction here,
which facilitates  understanding the proof of their
sharpness which we give in Section~\ref{sec:proof}.

For any integer $N$, define the monomial vector
\[
[x^N] = \bbm 1 & x_1 &\cdots & x_n & x_1^2 & x_1x_2 &\cdots &
x_n^N  \ebm^T.
\]
Then $[x^N][x^N]^T$ is a square matrix and we write
\[
[x^N][x^N]^T  = \sum_{0\leq |\af|\leq 2N} A_\af x^\af
\]
for some symmetric $0/1$-matrices $A_\af$. When $n=1$, the $A_\af$ are
Hankel matrices, and when $n>1$, the $A_\af$ are generalized
Hankel matrices. Suppose $\mu$ is a  nonnegative measure on
$\re^n$ with  total mass equal to one. Integrating the above
identity gives us
\[
M_N(y)= \int_{\re^n}  [x^N][x^N]^T  d \mu(x) = \sum_{0\leq
|\af|\leq 2N} A_\af y_\af
\]
where $y_\af = \int_{\re^n} x^\af d \mu(x)$ are the {\it moments
of $\mu$}. The matrix $M_N(y)$ is also called the {\it moment
matrix of order $N$}.

\subsection{SDP representation I}
\label{sec:SDPrepI}

Now we give the first construction of the lifted LMI which only
uses finitely many moments. Let $\mu(\cdot)$ be any nonnegative
measure such that $\mu(\re^n)=1$.
For any $\nu\in\{0,1\}^m$, define new polynomials
$g^\nu(x) := g_1^{\nu_1}(x) \cdots g_m^{\nu_m}(x)$. Let $d_\nu =
\lceil \deg(g_1^{\nu_1}\cdots g_m^{\nu_m})/2 \rceil$. For an fixed
integer $N \geq d_\nu$, define the {\it localizing moment matrix}
$M_{N - d_\nu}(g^\nu y)$  by
\[
M_{N - d_\nu}(g^\nu y) = \int_{\re^n}  g^\nu(x) [x^{N -
d_\nu}][x^{N - d_\nu}]^T \, d \mu(x) =
\sum_{0 \leq |\af| \leq 2N } A_\af^\nu y_\af
\]
where $y_\af = \int_{\re^n} x^\af d \mu(x)$ are the moments and
symmetric $A_\af^\nu$ are the coefficient matrices such that
\[
g^\nu(x) [x^{N - d_\nu}][x^{N - d_\nu}]^T\,  =  \, \sum_{0 \leq
|\af| \leq 2N } A_\af^\nu x^\af.
\]
For any integer $N \geq \max_\nu d_\nu = \lceil \deg(g_1 \cdots
g_m)/2\rceil$, if $\supp(\mu) \subseteq S$, then
\begin{align} \label{mLMI}
\forall \, \nu \in \{0,1\}^m, \quad M_{N-d_\nu}(g^\nu y) \succeq
0, \quad y_0 = 1.
\end{align}
Let $e_i$  denote the $i$-th unit vector in $\re^n$ whose only
nonzero entry is one and occurs at index $i$. If we set $y_0=1$
and $y_{e_i}=x_i$ in \reff{mLMI}, then it becomes the LMI
\begin{align} \label{LmiN}
\forall \, \nu \in \{0,1\}^m, \quad A_0^\nu+\underset{1\leq i\leq
n}{\sum} A_{e_i}^\nu x_i + \underset{1< |\af|\leq 2N}{\sum}
A_\af^\nu y_\af \succeq 0.
\end{align}
We mention that the LMI~\reff{LmiN} is essential the same as
the LMI~(2.11) in \cite{Las06}.

Let $\hS_N$ denote the set of all vectors $(x,y)$ satisfying
\reff{LmiN}. Note that $\hS_N \subset \re^n\times
\re^{^{\binom{n+2N}{n}-n-1} }$. For each $N$, define the
projection mapping
\[
\baray{rl}
\rho_N:\, \re^n \times \re^{^{\binom{n+2N}{n}-n-1} } & \rightarrow \re^n \\
(x,y) & \mapsto x. \earay
\]
The set $S$ is contained in the projection $\rho_N(\hS_N)$ of
$\hS_N$ onto $x$ space, because for any $x\in S$ the vector
$y=(y_\af)$ given by $y_\af = x^\af$ makes $(x,y)$ satisfies the
LMI~\reff{LmiN}. And obviously, the bigger $N$ is, the smaller the
projection $\rho_N(\hS_N)$ is. So, for any $N\geq \max_\nu d_\nu$,
we have the following chain
\[
\rho_N(\hS_N) \supseteq \rho_{N+1}(\hS_{N+1}) \supseteq \cdots
\supseteq  S.
\]
A natural question is whether there exists a finite integer $N$
such that $\rho_N(\hS_N) = S$.

\medskip
One typical approach to  this question is to use  linear
functionals to separate points in $\hat S_N $ from the convex set
$S$. Specifically, given a unit length vector $\ell \in \re^n$,
let $\ell^*$ be the minimum value of the linear function $\ell^Tx$
over the set $S$, let $u \in S$ denote the minimizer, which must
exist and be on the boundary $\pt S$. Since $S$ has nonempty
interior, the {\it Slater's condition} holds, and hence the {\it
first order optimality condition} is satisfied. So there exist
Lagrange multipliers $\lmd_1\geq 0, \cdots, \lmd_m\geq 0$ such
that $\ell =  \sum_i \lmd_i \nabla g_i(u)$. Suppose each
defining polynomial $g_i(x)$ is concave on $S$. Then \be
f_\ell(x):= \ell^Tx-\ell^*-\sum_i \lmd_i g_i(x) \ee is a convex
function such that $f_\ell(u)=0$ and $\nabla f_\ell(u) =0$. Thus
\[
  f_\ell(x) \geq f_\ell(u) + \nabla f_\ell(u)^T (x-u) =0,
  \quad \forall \, x\in S.
\]
In other words, $f_\ell$  is nonnegative on $S$ and so we could
wish $f_\ell(x)$ to have Schm\"{u}dgen's representation
\[
f_\ell(x) = \sum_{\nu\in\{0,1\}^m} \sig_\nu(x) g^\nu(x)
\]
for some particular SOS polynomials $\sig_\nu(x)$. Notice that
this representation is not implied by Schm\"{u}dgen's
Positivstellensatz \cite{Smg} because $f_\ell(x)$ has a zero point
$u$ on $S$.

Indeed, validating the lifted LMI $\hS_N$ for some finite integer
$N$ amounts to proving that for all $\ell$ the polynomial
$f_\ell(x)$ has Schm\"{u}dgen's representation with uniform
(in $\ell$) degree bounds on SOS polynomials $\sig_\nu(x)$; this we
will see in Section~\ref{sec:proof}.
This is equivalent to proving that a property on $S$ called
{\it Schm\"{u}dgen's Bounded Degree Representation (S-BDR)} of affine polynomials holds for $S$,
which means that there exists $N>0$ such that
for {\it almost every} pair $(a, b)\in \re^n \times \re$
\[
a^Tx +b > 0 \mbox{ on } S  \qquad \Rightarrow \qquad
a^Tx + b = \sum_{\nu\in\{0,1\}^m} \sig_\nu(x) g^\nu(x)
\]
for some SOS polynomials $\sig_\nu(x)$ with degree bounds
$\deg(\sig_\nu) + \deg(g^\nu) \leq 2N$.
When $S$ is a compact set (not necessarily convex),
Lasserre \cite{Las06} showed that if the S-BDR property holds
then the convex hull $\mbox{conv}(S)$ of $S$ equals
$\rho_N(\hS_N)$ for $N$ big enough.

S-BDR is a very nice restatement that the lift
\eqref{mLMI} and \eqref{LmiN} produces an SDP representation.
It reduces the main issue to finding
concrete and practical conditions assuring
the exactness of the lifted LMI~\reff{LmiN},
which is what we do in this paper.
Actually, when every polynomial $g_i(x)$ is concave on $S$
and strictly concave on $\bdS_i \cap \bdS$,
we can prove a stronger property called
{\it Schm\"{u}dgen's Bounded Degree
Nonnegative Representation (S-BDNR)}
of affine polynomials
holds for $S$, which means that there exists $N>0$ such that
for {\it every} pair $(a, b)\in \re^n \times \re$
\[
a^Tx +b \geq 0 \mbox{ on } S  \qquad \Rightarrow \qquad
a^Tx + b = \sum_{\nu\in\{0,1\}^m} \sig_\nu(x) g^\nu(x)
\]
for some SOS polynomials $\sig_\nu(x)$ with degree bounds
$\deg(\sig_\nu) + \deg(g^\nu) \leq 2N$.
As we can see, S-BDNR is a stronger property than S-BDR.
When S-BDNR property holds, S-BDR also holds
and hence $ S=\rho_N(\hS_N)$ for $N$
big enough by Theorem~2 in Lasserre \cite{Las06}.
To illustrate this we now state our
theorem for concave functions and S-BDNR,
while its proof will not  be given until
Section~\ref{sec:proofConvcaveg}.

\begin{theorem} \label{SeqLmiPrj}
Suppose $S=\{x\in\re^n:\, g_1(x)\geq 0, \cdots, g_m(x)\geq 0\}$ is
compact convex and has nonempty interior. Assume the $g_i(x)$ are concave
on $S$. For each $i$, if either $-\nabla^2 g_i(x)$ is SOS or
$-\nabla^2 g_i(u) \succ 0 $ for all $u \in \bdS_i \cap \bdS$, then
the S-BDNR property holds for $S$ and
there exists $N>0$ such that $ S=\rho_N(\hS_N)$.
\end{theorem}

\subsection{SDP representation II}
\label{sec:SDPrepII}

In LMI~\reff{LmiN}, the size of LMI is unfortunately exponential
in $m$. It is huge when $m$ is big. This is because we have used
all the possible products $g^\nu(x) = g_1^{\nu_1}(x)\cdots
g_m^{\nu_m}(x)$ for all index vector $\nu \in \{0,1\}^m$. If we
use only linear products, we can get a similar LMI
\begin{align} \label{putLMI}
\forall \, 0 \leq k \leq m, \quad A_0^{(k)}+\underset{1\leq i\leq
n}{\sum} A_{e_i}^{(k)} x_i + \underset{1< |\af|\leq 2N}{\sum}
A_\af^{(k)} y_\af \succeq 0
\end{align}
where symmetric matrices $A_\af^{(k)} = A_\af^{e_k}$ in
LMI~\reff{LmiN} ($e_0$ is the zero index vector).
We mention that the LMI~\reff{LmiN} is essential the same as
the LMI~(2.12) in \cite{Las06}.

Similar to
LMI~\reff{LmiN}, let $\tS_N$ be the set of all vectors $(x,y)$
satisfying \reff{putLMI} and $\tilde \rho_N$ be the projection
mapping into $x$-space. Then the following chain relation again
holds
\[
\tilde \rho_N(\tS_N) \supseteq \tilde \rho_{N+1}(\tS_{N+1})
\supseteq \cdots \supseteq  S.
\]
The natural question is whether $\tilde \rho_N(\tS_N)=S$ for some
finite integer $N$. This can  be shown true under the so called
{\it archimedean condition}: There exist SOS polynomials
$s_0(x),s_1(x), \cdots, s_m(x)$ and a number $R>0$ big enough such
that
\[
R- \sum_{i=1}^n x_i^2  =  s_0(x)+s_1(x)g_1(x)+ \cdots +
s_m(x)g_m(x).
\]
Note that the archimedean condition implies $S$ is compact. But
the converse might not be true. However, a compact $S$ can be
forced to satisfy the archimedean condition by adding a
``redundant" constraint like $R-\sum_{i=1}^m x_i^2 \geq 0$ for
sufficiently large $R$.

Similar to the lifted LMI~\reff{LmiN},
validating the exactness of the lifted LMI $\hS_N$
amounts to proving that for every $\ell$ the polynomial $f_\ell(x)$
has Putinar's representation
\[
f_\ell(x) = \sig_0(x) + \sig_1(x) g_1(x) + \cdots + \sig_m(x) g_m(x)
\]
with uniform degree bounds on SOS polynomials $\sig_i(x)$.
This is equivalent to proving that the so-called
{\it Putinar-Prestel's Bounded Degree Representation
(PP-BDR)} property \cite{Las06}
holds for $S$, that is, there exists $N>0$ such that
for {\it almost every} $(a,b)\in \re^n \times \re$
\[
a^Tx+b > 0 \mbox{ on } S  \qquad \Rightarrow \qquad
a^Tx + b = \sig_0(x) + \sig_1(x) g_1(x) + \cdots + \sig_m(x) g_m(x)
\]
for some SOS polynomials $\sig_i(x)$ with degrees
$\deg(\sig_i) + \deg(g^\nu) \leq 2N$.
When $S$ is a compact set (not necessarily convex),
Lasserre \cite{Las06} showed that if the PP-BDR property holds
then the convex hull $\mbox{conv}(S)$ of $S$ equals
$\tilde \rho_N(\tS_N)$ for $N$ big enough.

As with S-BDR, the PP-BDR property is a nice reformulation
of the exactness of the lifted LMI~\reff{putLMI} and
is usually
not directly checkable.
In this paper,
under the archimedean condition,
when every polynomial $g_i(x)$ is concave on $S$
and strictly concave on $\bdS_i \cap \bdS$,
we can prove a property called
{\it Putinar-Prestel's Bounded Degree Nonnegative Representation (PP-BDNR)}
of affine polynomials
holds for $S$, which means that there exists $N>0$ such that
for {\it every} pair $(a, b)\in \re^n \times \re$
\[
a^Tx+b \geq 0 \mbox{ on } S  \qquad \Rightarrow \qquad
a^Tx + b = \sig_0(x) + \sig_1(x) g_1(x) + \cdots + \sig_m(x) g_m(x)
\]
for some SOS polynomials $\sig_i(x)$ with degrees
$\deg(\sig_i g_i) \leq 2N$ (we denote $g_0(x)=1$).
Obviously, PP-BDNR implies PP-BDR.
When PP-BDNR property holds, PP-BDR also holds
and hence $ S=\tilde \rho_N(\tS_N)$ for $N$
big enough by Theorem~2 in Lasserre \cite{Las06}.
The following theorem,
whose proof will be given in Section~\ref{sec:proofConvcaveg},
illustrates the this for concave functions and PP-BDNR.

\begin{theorem} \label{PutLmiLift}
Suppose $S=\{x\in\re^n:\, g_1(x)\geq 0, \cdots, g_m(x)\geq 0\}$ is
compact and has nonempty interior. Assume the $g_i(x)$ are concave
on $S$ and the archimedean condition holds. For each $i$, if
either $-\nabla^2 g_i(x)$ is SOS or $-\nabla^2 g_i(u) \succ 0 $ for
all $u \in \bdS_i \cap \bdS$, then
the PP-BDNR property holds
and there exists $N>0$ such that $
S=\tilde \rho_N(\tS_N)$.
\end{theorem}

\medskip

It is possible that the defining polynomials $g_i(x)$ are not
concave but the set $S$ is still convex.
In this case, does $S$
have an SDP representation? After some modifications in
LMI~\reff{LmiN}, the answer is affirmative in very general
situations, which is our Theorem~\ref{thm:mainRep}.
However, our
proof of Theorem \ref{thm:mainRep} uses Theorem \ref{SeqLmiPrj} or
\ref{PutLmiLift} as a stepping stone.

\section{The SDP representation when $g_i(x)$ are sos-concave}
\label{sec:globconcav} \setcounter{equation}{0}

Lasserre \cite{Las06}  and Parrilo \cite{Par06}
proposed  recipes for an SDP representation. In this section we
give a sufficient condition such that the LMI constructed in
Lasserre \cite{Las06} is a lifted LMI of $S$.
We assume the polynomials $g_i(x)$ are concave (not necessarily
strictly concave) in the whole space $\re^n$. Certainly the set
$S=\{x\in \re^n:\, g_1(x)\geq 0, \cdots, g_m(x)\geq 0\}$ is
convex.

\medskip
As was shown in Lasserre \cite{Las06}, the set $S$ is contained in
the projection of $\hS$ defined by LMI
\begin{align}
\label{lassLMI} \left. \baray{r}
M_{d_g}(y) \succeq 0 \\
L_{g_1}(y),\cdots, L_{g_m}(y)\geq 0\\
y_0 = 1 \earay\right\}
\end{align}
where $d_g=\max_i \lceil \deg(g_i)/2\rceil$ and $
L_{g_i}(y)=\sum_\af g_\af^{(i)}y_\af$ if we write $g_i(x)=\sum_\af
g_\af^{(i)}x^\af$. The projection onto $x$ -space of $\hS$ is $\{
x:  \, \exists  \  y \in \hS, \  \  x_i = y_{e_i}, \, 1\leq i\leq
n \}$. It is natural to ask whether the LMI~\reff{lassLMI} is a
lifted LMI for $S$, i.e., the projection of LMI~\reff{lassLMI}
onto $x$-space equals $S$.

The standard approach to  this question is to use separating
linear functionals. As we did in Section~\ref{sec:LMIlift}, for
each vector $\ell \in \re^n$, let $\ell^*$ be the minimum value of
$\ell^Tx$ over the set $S$, $u \in S$ be the minimizer, which must
be on the boundary $\pt S$. If there is some point in  the
interior of $S$, then the Slater's condition holds and hence there
exist Lagrange multipliers $\lmd_1\geq 0, \cdots, \lmd_m \geq 0$
such that the optimality condition $\ell = \sum_i \lmd_i \nabla
g_i(u)$ holds, and hence \be f_\ell(x):= \ell^Tx-\ell^*-\sum_i
\lmd_i g_i(x) \ee is a convex and nonnegative polynomial in the
whole space $\re^n$ such that $f_\ell(u)=0$ and $\nabla
f_\ell(u)=0$ (see \cite{Las06}). Under the assumption that the
polynomial $f_\ell$ is SOS for {\it every} vector $\ell\in\re^n$,
Lasserre \cite{Las06} showed the LMI~\reff{lassLMI} is a lifted
LMI for $S$. If $f_\ell$ is not a sum of squares for some
particular $\ell\in\re^n$, then the LMI~\reff{lassLMI} might not
be a lifted LMI for $S$.

\medskip

Although it is very difficult to tell if a polynomial is
nonnegative, it is more tractable to check if a polynomial is SOS,
which can be done by solving an SDP feasibility problem, e.g., by
softwares like {\it SOSTOOLS} \cite{PPSP} and {\it Gloptipoly}
\cite{Gloptipoly}. However, it is impossible to check if $f_\ell$
is SOS for uncountably many vectors $\ell\in\re^n$.
Here we give a sufficient condition for the LMI~\reff{lassLMI} to
be a lifted LMI of $S$,
which can be checked numerically.
Let us start the discussion with a lemma.

\medskip

\begin{lem} \label{lem:SosInt}
If a symmetric matrix polynomial $P(x) \in \re[x]^{r\times r}$ is
SOS, i.e., $P(x)=W(x)^TW(x)$
for some possibly nonsquare matrix polynomial
$W(x)\in \re[x]^{k\times r}$,
then for any $u\in \re^n$ the double integral
\[
\int_0^1 \int_0^t  P(u+s(x-u)) \, ds\, dt
\]
is also a symmetric SOS matrix polynomial in $\re[x]^{r\times r}$.
In particular, when $r=1$, the above double integral of scalar SOS
polynomials is also SOS.
\end{lem}
\noindent {\it Proof} \,
Let $2d = \deg(P)$. Let $\xi \in \re^r$ be a symbolic vector. Then
$P(x)$ is SOS if and only if
\[
\xi^TP(x)\xi = [\xi x^d]^T A^T A [\xi x^d]
\]
for matrix $A$. Here $[\xi x^d]$ denotes the vector of monomials
\[
\big[ \,\xi_1 [x^d]^T \quad \cdots  \quad \xi_r [x^d]^T \, \big]^T.
\]
Note that $\xi^TP(x)\xi$ has degree $2$ in $\xi$. If $P(x)$ is SOS,
we can assume the above $A$ exists. In monomial vector
$[\xi x^d]$ we replace $x$ by $u+s(x-u)$. Each entry of
$[\xi (u+s(x-d))^d]$ is a polynomial in $x$ whose coefficients are
polynomials in $u$ and $s$. So there exists a matrix polynomial
$C(u,s)$ such that $[\xi (u+s(x-d))^d] = C(u,s) [\xi x^d] $. Therefore we
have
\[
\int_0^1 \int_0^t  \xi^TP(u+s(x-u))\xi \, ds\, dt = [\xi x^d]^T  B^TB [\xi x^d]
\]
where $B$ is a matrix such that
\[
\int_0^1 \int_0^t  C(u,s)^TA^TAC(u,s)  \, ds\, dt = B^TB.
\]
Therefore, the double integral of matrix polynomial in the lemma is SOS.
\eproof


\medskip
\noindent
{\it Remark:} The above integral is the limit of a
sequence of Riemann sums, which are all SOS with bounded degrees.
So intuitively the integral must also be SOS.

\medskip

\begin{lem} \label{lem:LPhess}
Let $p(x)$ be polynomial such that $p(u)=0$ and $\nabla p(u)=0$
for some point $u \in \re^n$. If the Hessian $\nabla^2p(x)$ is
SOS, then $p(x)$ is SOS.
\end{lem}
\noindent {\it Proof} \,
Let $q(t)=p(u+t(x-u))$ be a univariate polynomial in $t$. Then
\[
q^{\prm\prm}(t) = (x-u)^T \nabla^2 p(u+t(x-u))(x-u).
\]
So we have
\[
p(x) = q(1) = (x-u)^T \left( \int_0^1 \int_0^t \nabla^2
p(u+s(x-u)) \,  ds \, dt \, \right) (x-u).
\]
Since $\nabla^2 p(x)$ is SOS, the middle double integral above
should also be SOS, by Lemma~\ref{lem:SosInt}.
So $p(x)$ is also SOS.
\eproof

\medskip

The following theorem gives the sufficient condition which we are
aiming at.

\begin{theorem}
Assume $S$ has nonempty interior. If every $g_i(x)$
is sos-concave, then \reff{lassLMI} is a lifted LMI representation
for $S$.
\end{theorem}
\noindent {\it Proof} \,
It is obviously that $S$ is contained
in the projection of the set $\hat S$.
If they are not equal,
there must exist some $\hat y \in \hat S$ such that
$\hat x = (\hat y_{e_1}, \cdots,  \hat y_{e_n}  ) \notin S$.
Since $S$ is closed, there exists a supporting hyperplane
of $S$ that excludes $\hat x$, i.e.,
there exists a unit length vector $\ell \in \re^n$ such that
\[
\ell^Tx -\ell^* \geq 0,  \,\,\, \forall \, x \in S,
\quad \ell^T u -\ell^* = 0, \,\,\, \exists \, u \in \bdS,
\qquad \ell^T \hat x - \ell^* < 0.
\]
Then $u$ is a minimizer for the above.
Since $S$ has nonempty interior, the Slater's condition holds, and hence
there must exist Lagrange multipliers $\lmd_1\geq 0, \cdots,
\lmd_m\geq 0$ such that $\ell = \sum_{i=1}^m \lmd_i \nabla
g_i(u)$. Thus $f_\ell(x) = \ell^Tx-\ell^*-\sum_i \lmd_ig_i(x)$ is
convex polynomial in $\re^n$ such that $f_\ell(u)=0$ and $\nabla
f_\ell(u)=0$.
Note that $\nabla^2 f_\ell(x) = \sum_{i=1}^m
 \lmd_i (-\nabla^2 g_i(x))$ is SOS,
since all $g_i(x)$ are sos-concave.
Hence Lemma~\ref{lem:LPhess}
implies that $f_\ell(x)$ is SOS,
i.e., there exists a symmetric matrix $W\succeq 0$ such that the identity
\[
\ell^Tx-\ell^* = \overset{m}{\underset{i=1}{\sum}} \lmd_i g_i(x) +  [x^{d_g}]^T W [x^{d_g}]
\]
holds. In the above identity, replace each monomial $x^\af$ by $\hat y_\af$,
then we get
\[
\ell^T \hat x-\ell^* = \overset{m}{\underset{i=1}{\sum}}
 \lmd_i L_{g_i}(\hat y)
+ Trace \big(W \cdot M_{d_g}(\hat y) \big) \geq 0,
\]
which contradicts the previous assertion $\ell^T \hat x - \ell^* < 0$.
\eproof

\noindent
{\it Remarks:}
(i) We do not need assume $S$ is compact whenever
every defining polynomial $g_i(x)$ is assumed to be sos-concave.
(ii) Checking whether $g_i(x)$ is sos-concave can be done numerically.
Obviously, $-\nabla^2g_i(x)$ is SOS if and only if the polynomial
$ -\sum_{k,\ell=1}^n \frac{\pt^2 g_i(x)}{\pt x_k\pt x_\ell} \xi_k
\xi_\ell$ in $(x, \xi)$ is SOS. This can be checked numerically by
solving an SDP feasibility problem, e.g., by softwares {\it
SOSTOOLS} and {\it Gloptipoly}.

\section{Concave defining functions for poscurv-convex sets }
\label{sec:ConDef} \setcounter{equation}{0}

It is possible that the set $S_i=\{x\in\re^n:\, g_i(x)\geq 0\}$ is
convex but that its defining polynomials
$g_i(x)$ are neither concave nor quasi-concave.
The goal of this section is to find
a new set of defining polynomials for $S$.
Under some modest hypotheses,
we show that $S$ can be defined by a set of
polynomials having negative definite Hessians on $S$
when $S$ is defined by strictly quasi-concave polynomials,
and that there is a
new set of (local) defining polynomials $p_i$ for the $S_i$ which
are strictly concave.
These results, Proposition~\ref{prop:nicefunc},
Proposition~\ref{prop:approxp}
and Corollary \ref{cor:approxp},
might be of interest independent of our SDP representation
application.
Proposition~\ref{prop:nicefunc} and Proposition~\ref{prop:approxp}
allow us to prove Theorems~\ref{thm:mainRep} and \ref{thm:posCurv}
by using the $p_i$ together with the given $g_i$.

We start with a short review of curvature and convexity.
For each
$u\in Z(g_i)$ with $\nabla g_i(u) \not = 0$, the hyperplane
\[
H_i(u)=\{x\in \re^n: \nabla g_i(u)^T(x-u)=0\}
\]
is tangent to $Z(g_i)$ at $u$. Also the set $S_i$ is strictly convex,
in some neighborhood of  $u$ if and only if
\be
\label{eq:stcon}
H_i(u) \cap  Z(g_i) \cap B(u,\eps)  = \{ u\}
\ee
for a small $\eps >0$,
which holds if and only if
\begin{align*}
g_i(x) =  (x-u)^T \nabla^2 g_i(u) (x-u) + o(\|x-u\|^2) \ < 0,
\qquad
  \forall \, u\ne x \in \, H_i(u)\cap B(u,\eps).
\end{align*}
This implies the quadratic form associated with the negative
Hessian
$$\Phi(v):=- v^T\nabla^2 g_i(u)v $$
is nonnegative for all $v$ in the  tangent space
\[
\nabla g_i(u)^\perp = \{ v \in \re^n:\, \nabla g_i(u)^Tv = 0 \}.
\]

We follow common usage (c.f.
\cite{EOM} )and call the quadratic function $\Phi$ on the tangent
space $\nabla g_i(u)^\perp$  the
{\it second fundamental form} of $Z(g_i)$.
A surface has {\it  positive curvature} at $u$ if
and only if the second fundamental form is strictly positive
definite there.
Obviously, when $S_i$ is convex
 positive curvature of $Z(g_i)$ everywhere on $\bdS_i$ implies strict
convexity of $S_i$, but the converse is not necessarily true.
Results like this which assume nonnegative curvature,
but do not require $S_i$ to be convex as a
hypothesis are in \cite{Dmitriev}.

Note that while
in this paper we sometimes require
$\nabla g_i(u)\not = 0$  at $u\in Z(g_i) \cap S$,
the definition of positive curvature of $Z(g_i)$
itself does not.
Indeed if the gradient $\nabla g_i(u)$ vanishes,
then one can interpret
$\nabla g_i(u)^\perp$ as the whole space and
its
negative Hessian $-\nabla^2 g_i(u)$
is required to be positive definite.

We will distinguish the cases that $S$ is defined
by strictly quasi-concave functions
and $S$ has positively curved boundary.
A new set of defining polynomials for $S$
will be discussed in two subsections separately.

\subsection{ Convex sets defined by quasi-concave functions }

The following proposition gives $S_i$ a new defining polynomial
$p_i$
whose Hessian is negative definite on $S$
when $g_i(x)$ is strictly quasi-concave on $S$.

\begin{prop} \label{prop:nicefunc}
Assume $g_i(x)$ is strictly quasi-concave on $S$. Then there
exists a polynomial $h_i(x)$ such that $h_i(x)$ is positive on $S$
and the product $p_i(x) = g_i(x)h_i(x)$ is concave on $S$ and has
negative definite Hessian there.
\end{prop}

We give the proof of this proposition after introducing some
lemmas. Without loss of generality, suppose $ 0 \leq g_i(x) \leq
1$ for all $x\in S$, since $S$ is compact, because otherwise we
can scale the coefficients of $g_i(x)$. The set $S_i$ can be
convex without $g_i(x)$ being concave. So the Hessian $-\nabla^2
g_i(x)$ might be indefinite. However, the Hessian $-\nabla^2
g_i(x)$ can have at most one negative eigenvalue for $x\in S$, and
better yet the Hessian can be ``relaxed" to yield the ``modified
Hessian'' which is positive definite.

\begin{lem} \label{modHes}
Assume $g_i(x)$ is strictly quasi-concave on $S$. Then we have
\bnum \item [(a)] There exists $M$ sufficiently large such that
the {\bf modified Hessian}
\[ - \nabla^2 g_i(x) + M \nabla g_i(x) \nabla g_i(x)^T \]
is positive definite for all $x\in S$.

\item [(b)] If $g_i(x)$ is concave, then the above modified
Hessian is positive definite for any $M >0$. \enum
\end{lem}
\noindent {\it Proof} \,
(a). Let $U_i = \{ x\in S:\, -\nabla^2 g_i(x) \succ 0\}$, which is
an open set in $S$. Then $V_i= S-U_i$ is a compact set. Choose an
arbitrary point $u\in V_i$ and let $\af = g_i(u) \geq 0$. Then
$Z(g_i-\af) = \{ x\in \re^n:\, g_i(x) = \af\}$ has positive
curvature in $S$.
Note that for every $u\in V_i$, the negative Hessian $-\nabla^2
g_i(u)$ is not positive definite and hence $\nabla g_i(u)\ne 0$,
because otherwise $\nabla g_i(u)^\perp$ is the whole space $\re^n$
which implies $-\nabla^2 g_i(u)$ is positive definite.

Let $Q_i(u) = [\, \nabla g_i(u) \,\,\,\, \hat Q_i(u) \,] \in
\re^{n\times n}$ be a nonsingular matrix such that $ \nabla g_i(u)^T \hat
Q_i(u) =0$. Then
\begin{align*}
& \quad Q_i(u)^T \left(- \nabla^2 g_i(u) + M \nabla g_i(u) \nabla
g_i(u)^T\right) Q_i(u) \\
= & \bbm - \nabla g_i(u)^T \nabla^2 g_i(u) \nabla g_i(u) + M
\|\nabla g_i(u)\|^2
&  - \nabla g_i(u)^T \nabla^2 g_i(x) \hat Q_i(u) \\
- \hat Q_i(u)^T \nabla^2 g_i(u)\nabla g_i(u) & - \hat Q_i(u)^T
\nabla^2 g_i(u)\hat Q_i(u) \ebm.
\end{align*}
For $u\in V_i$,
 $- \hat Q_i(u)^T \nabla^2 g_i(u)\hat Q_i(u) \succ 0$ and
$\nabla g_i(u) \ne 0$ on $V_i$. Since $V_i$ is compact, we can
choose $M$ big enough such that the modified Hessian is positive
definite for all $u\in V_i$. When $u\in U_i$, the modified Hessian
is obviously positive definite.

\medskip

(b). If $g_i(x)$ is concave, then the modified Hessian is
obviously positive semidefinite. We need show it is positive
definite for any $ M >0$. Suppose for some $u\in S$ and a vector
$\xi \in \re^n$
\[
- \xi^T\nabla^2 g_i(u)\xi + M \xi^T\nabla g_i(u) \nabla
g_i(u)^T\xi =0.
\]
Then it must hold
\[
-\xi^T\nabla^2 g_i(u)\xi =0, \qquad  \nabla g_i(u)^T\xi =0.
\]
Since $-\nabla^2 g_i(u) \succ 0$ in the tangent space $\nabla
g_i(u)^\perp$, we must have $\xi=0$, which completes the proof.
\eproof

\begin{exm}
(1) The following set is strictly convex
\[
\{x\in \re^2:\, \underbrace{x_1x_2 -1 }_{g_1(x)} \geq 0, \,
\underbrace{1-(x_1-1)^2-(x_2 -1)^2 }_{g_2(x)} \geq 0 \}.
\]
$g_2(x)$ is strictly concave, but $g_1(x)$ is not concave.
However, for any $ M> \half$, the modified Hessian
\[-\nabla^2 g_1(x) + M \nabla g_1(x) \nabla g_1(x)^T \]
is positive definite on $S$. \\  \vskip .01cm \noindent (2) The
condition that $g_i$ is strictly quasi-concave in $S$ in
Lemma~\ref{modHes} can not be weakened to $S$ is strictly convex.
For a counterexample, consider the strictly convex set
\[
\{x\in \re^2:\, \underbrace{x_2 -x_1^3 }_{g_1(x)} \geq 0, \,
\underbrace{x_2+x_1^3 }_{g_2(x)} \geq 0, \,
\underbrace{1-(x_1-1)^2-(x_2 -1)^2 }_{g_3(x)} \geq 0 \}.
\]
For $i=1,2$, no matter how big $M$ is, the modified Hessian
\[
-\nabla^2 g_i(x) + M \nabla g_i(x) \nabla g_i(x)^T = \bbm  \pm
6x_1 + 9 M x_1^4 &  \pm 3 x_1^2 \\ \pm 3x_1^2 & M  \ebm
\]
can not be positive semidefinite near the origin.
\end{exm}

\begin{lem} \label{lem:phi}
For an arbitrarily large number $M>0$, there exists a univariate
polynomial function $\phi(t)$ such that for all $t\in[0,1]$
\begin{align} \label{mulpos}
\phi(t)>0, \qquad \phi(t) + \phi^{\prm}(t) t > 0, \qquad \frac{ 2
\phi^\prm(t) + \phi^{\prm\prm}(t) t } {\phi(t) + \phi^{\prm}(t) t}
\leq -M .
\end{align}
\end{lem}
\noindent {\it Proof} \,
The smooth function $\psi(t) = \frac{1-e^{-(M+1)t}}{(M+1)t}$
satisfies the following
\begin{align*}
\psi(t) + t \psi^{\prm}(t) &= (t \psi(t))^{\prm} = e^{-(M+1)t} \\
2\psi^{\prm}(t) + t \psi^{\prm\prm}(t) & = (\psi(t) + t
\psi^{\prm}(t) )^{\prm} = -(M+1) e^{-(M+1)t}.
\end{align*}
So $\psi(t)$ satisfies \reff{mulpos}. Let $\psi(t) =
\sum_{k=0}^\infty a_k t^k$ be the power series expansion, and let
$\psi_N(x) = \sum_{k=0}^N a_k t^k$ be the truncated summation.
Note that $\psi_N$ converges to $\psi$ uniformly on $[0,1]$. For
arbitrarily small $\vareps >0$, we can choose $N$ big enough such
that for all $t\in [0,1]$
\[
|\psi_N(t)-\psi(t)| <\vareps, \qquad
|\psi_N^{\prm}(t)-\psi^{\prm}(t)| <\vareps, \qquad
|\psi_N^{\prm\prm}(t)-\psi^{\prm\prm}(t)| <\vareps.
\]
Then the polynomial $\phi(t) = \psi_N(t)$ satisfies \reff{mulpos}
when $N$ is big enough.
\eproof

\medskip

\noindent {\it Proof of Proposition \ref{prop:nicefunc} } \,
Let $\phi(t)$ be a polynomial satisfying \reff{mulpos} and $h_i(x) = \phi(g_i(x))$,
which is positive on $S$, since $g_i(S)\subseteq [0,1]$. Then a
direct calculation shows for $p_i(x)=g_i(x)h_i(x)$
\begin{align*}
-\nabla^2( p_i(x)) & = - \left(\phi(g_i) + \phi^{\prm}(g_i) g_i
\right)\nabla^2 g_i(x) + \left(2 \phi^\prm( g_i ) +
\phi^{\prm\prm}( g_i ) g_i \right)
\nabla g_i(x) \nabla g_i(x)^T  \\
& = \left( \phi(g_i) + \phi^{\prm}(g_i) g_i  \right) \left(
-\nabla^2 g_i(x) - \frac{2 \phi^\prm( g_i ) + \phi^{\prm\prm}( g_i
) g_i}{\phi(g_i) + \phi^{\prm}(g_i) g_i} \nabla g_i(x) \nabla
g_i(x)^T \right).
\end{align*}
If $M$ in \reff{mulpos} is chosen big enough, by
Lemma~\ref{modHes} and Lemma~\ref{lem:phi}, the negative Hessian
$-\nabla^2 (p_i(x))$ must be positive definite for all $x\in S$.
\eproof

\medskip
\noindent
{\it Remark:} From the proof,
we can see that
both Proposition~\reff{prop:nicefunc} and Lemma~4.2
remain true if $S$ is replaced by any compact set $\Omega$
which is not convex or even connected.

\subsection{Convex sets with boundary having positive curvature}

This subsection ultimately shows that such an extendable
poscurv-convex set has  a very well behaved
defining function.
Recall $S_i$ is {extendable poscurv-convex with respect to} $S$ if
$g_i(x) >0$ whenever $x\in S$ lies in  the interior of $S_i$
and there exists a
poscurv-convex set $T_i \supseteq S$ such that
$\pt T_i\cap S = \bdS_i \cap S$.
First we give a result which says a poscurv-convex set can be defined
by a strictly concave smooth function.

\medskip

\begin{prop}
\label{prop:ConDefFun}
Suppose $T_i$ is a poscurv-convex set with the origin in the interior.
Then there is a
function $G_i(x)\, \in \,C^2(\re^n) \cap C^\infty(\re^n-\{0\})$
such that $T_i=\{x\in\re^n:\, G_i(x) \geq 0\}$, $\pt
T_i=\{x\in\re^n:\, G_i(x) = 0\}$,
$\nabla G_i(x) \ne 0$ for all $x\in \pt T_i$,
and $G_i(x)$ has negative
definite Hessian on $T_i$.
\end{prop}

\noindent {\it Proof} \,
Our strategy is to build a concave defining function for $T_i$ and
then to approximate it by a concave smooth function. This takes
several steps.
Since $T_i$ is compact, for any $0\ne x\in \re^n$, there exists a
unique positive scalar $\af(x)$ such that
$\frac{1}{\af(x)}x=:r(x)$ lies
on the boundary $\pt T_i$. Define $\af(0)=0$. Indeed $\alpha $ is
the classical Minkowski function (\cite{C}), and $T_i
=\{x\in\re^n:\, \af(x) \leq 1 \} $. The function $\af(x)$ is
convex.
Note we can write $x = \af(x) r(x)$ and
$\af(x)$ is smooth at $x\ne 0$,
because the boundary $\pt T_i$ is smooth.

Let $\tG(x)=1-\af(x)^3$. Thus $\tG(x)$ is a concave function and
is smooth everywhere except at 0. Moreover, the super level sets
satisfy
$$\left\{ x : \ \tG(x)\geq c \right\}
=\left\{ x : \ 1 \geq
\alpha\left(\frac{x}{\sqrt[3]{1-c}}\right)\right\} = \left\{x: \
\frac{x}{\sqrt[3]{1-c}} \in T_i \right\}=  \sqrt[3]{1-c} T_i
$$
for all $0 \leq c < 1$. Since $\pt T_i$ has positive curvature,
$\sqrt[3]{1-c}\, \pt T_i$ also has positive curvature. In summary,
the function $\tG$ is concave, strictly quasi-concave and smooth
except at $x=0$. However, we need a function that is twice
continuously differentiable on $T_i$ and has negative definite
Hessian there. The following produces one.
\smallskip

\noindent {\bf Claim:}  {\it \ \ $G_i(x):=  (1 - \eps \| x \|^2
)[1 -( \af(x) )^3 ] \, \in \,C^2(\re^n) \cap C^\infty(\re^n-\{0\})
$ has negative definite Hessian on $T_i$
when $\eps$ is small enough.} \\

\noindent {\it Proof of the Claim}.
Let $\psi(t) := 1-t^3$ and
then $ \tG:= \psi \circ \af $.
So at $x\ne 0$
\be
\label{eq:phiGderiv} \nabla \tG  = \psi' (\af) \nabla \af, \qquad
\nabla^2 \tG = \psi''( \af )  \nabla \af \nabla \af^T + \psi'
(\af) \nabla^2 \af.
\ee
Note that $\nabla \tG(x) \ne 0$ for all $x\in \pt T_i$,
since $\pt T_i$ is smooth.
Now we use the above to prove at $x\neq 0$ the
Hessian $ \nabla^2\tG(x)$ is negative definite. Obviously for
$0\ne x\in T_i$
\[
 \nabla \af \nabla \af^T \succeq 0, \qquad \psi''(\af) < 0,\qquad
 \nabla^2 \af \succeq 0 ,\qquad  \psi' (\af)  < 0
\]
and \reff{eq:phiGderiv} has the form of the modified Hessian. Thus
part (b) of Lemma \ref{modHes} implies $ \nabla^2 \tG(x)$ is
negative definite at $0\ne x\in T_i$. From $x =\af(x) r(x)$, we
have
\[
|\af(x)| = \frac{ |x^T r(x)|}{\|r(x)\|^2} \leq  \frac{
\|x\|}{\|r(x)\|} .
\]
For $x\ne 0$, $r(x)$ is on the boundary $\pt T_i$ and hence
$\|r(x)\| \geq \dt$ for some constant $\dt >0$. Thus $\af(x) =
\mc{O}(\|x\|)$ and then $\af(x)^3 = \mc{O}(\|x\|^3)$. So we can
see $\af(x)^3$ is at least twice differentiable at the origin; its
gradient and Hessian vanish there, and so do those of $\tG(x)$.
The function $\tG$ has negative definite Hessian except at $x=0$.
Obviously $G_i(x)\, \in \,C^2(\re^n) \cap C^\infty(\re^n-\{0\})$
and $\nabla G(x) \ne 0$ for all $x\in \pt T_i$.

To achieve strict concavity at $0$ take $ G_i(x): = (1 - \eps
\|x\|^2) \tG(x)$. Then
\[
\nabla^2 G_i(x) =  (1 - \eps \|x\|^2) \nabla^2 \tG(x) -2\eps
\left( \nabla \tG(x) x^T+ x\nabla \tG(x)^T\right) -2\eps \tG(x)
I_n.
\]
At $x=0$, $\nabla^2 G_i(0)=-2\eps$. Thus for $\eps >0$ the Hessian
of $G_i$ at $x=0$ is negative definite. We can take $\eps$ small
enough to keep the Hessian of $G_i$  negative definite on the
compact set $T_i$ away from $x=0$ as well, which completes the
proof of the claim.

\medskip

Obviously, $x\in T_i$ if and only if $G_i(x)\geq 0$, and $x\in \pt
T_i$ if and only if $G_i(x) =0$.
\eproof

\def\tx{\tilde x}

\begin{lem}
\label{lem:gG}
Assume $S=\{x\in\re^n:\, g_1(x) \geq 0, \cdots, g_m(x) \geq 0\}$
is compact convex and
$S_i=\{x\in\re^n:\, g_i(x) \geq 0\}$ is
extendable poscurv-convex with respect to $S$.
Then we have
\begin{enumerate}
\item [(i)] $\nabla g_i(x)$ does not vanish on the boundary
$\bdS_i \cap S$, and hence $\bdS_i \cap S$ is smooth. \item
[(ii)] Let $G_i(x)$ be the defining function for $T_i$ given by
Proposition~\ref{prop:ConDefFun}. Then
\[
w(x):= \frac{G_i(x)}{g_i(x)} \in C^2(S)\cap C^\infty(\bdS), \qquad
\text{ and } \qquad w(x) > 0 , \, \forall x \in S.
\]
\end{enumerate}
\end{lem}

\noindent {\it Proof} \,
(i) We show $\nabla g_i(x) \ne 0$ on the boundary $\bdS_i\cap S$.
For a contradiction, suppose $\nabla g_i(u)=0$ for some $u\in
\bdS_i \cap S$. Since $\bdS_i$ has positive curvature, we have
$-\nabla^2 g_i(u) \succ 0$. By continuity, $-\nabla^2 g_i(x) \succ
0$ when $x \in B(u,\eps)$ for some $\eps >0$. Since $S$ is convex
and has nonempty interior, there exists $v\in B(u,\eps)$ in the
interior of $S$. Thus
\[
g_i(v) < g_i(u) + \nabla g_i(u)^T (v-u) = g_i(u) =0.
\]
which contradicts $g_i(v) \geq 0$ since $v\in S$.
\medskip

(ii)
By assumption,
let $T_i\supseteq S$ be a poscurv-convex set
such that $\pt T_i\cap S = \bdS_i \cap S$;
thus  and $\pt T_i$ is nonsingular
and has positive curvature.
Without loss of generality, assume the
origin is in the interior of $S$. Then apply
Proposition~\ref{prop:ConDefFun} to $T_i$ and get a concave
defining function for $T_i$ such that $G_i(x)\in C^2(\re^n)\cap
C^\infty(\bdS)$ and it has a negative definite Hessian on $T_i
\supseteq S$. Similarly we can prove $\nabla G_i(x) \ne 0$ on the
boundary $\pt T_i$.

Now we need to show $w(x)= \frac{G_i(x)}{g_i(x)}$ is positive on
$S$ and belongs to $C^2(S)\cap C^\infty(\bdS)$.
Obviously it is
smooth in the interior or exterior
of $S_i$ except at $x=0$,
and twice differentiable at $x=0$. We need to show $w$ is smooth
on the boundary $\bdS_i \cap S$. Now fix a $u\in \bdS_i\cap S$.
Since $\nabla g_i(u) \ne 0$, we can find a local coordinate
transformation $x-u = t \nabla g_i(u) + B y$ to new coordinates
$(t, y)$ in $\re \times \re^{n-1}$. Here $B$ is a matrix such that
$\nabla g_i(u)^T B=0$ and $[\nabla g_i(u) \,\, B]$ is invertible.
The point $u$ corresponds to $(0,0)$ in the new coordinate. Then
apply the Taylor series expansion at point $u$ and get
$G_i(x)=G_i(t,y) = \sum_{k=1}^\infty a_k(y) t^k$ and
$g_i(x)=g_i(t,y) = \sum_{k=1}^\infty b_k(y) t^k$ for some smooth
scalar functions $a_k(y),b_k(y)$. The fact $\nabla g_i(u) \ne 0$
and $\nabla G_i(u) \ne 0$ implies
 $a_1(0) \ne 0$ and $b_1(0) \ne 0$.
Thus we can see
\[
w(x) = \frac{a_1(y) + \sum_{k=2}^\infty a_k(y) t^{k-1}} {b_1(y) +
\sum_{k=2}^\infty a_k(y) t^{k-1}}
\]
is smooth at $u$. Note that $a_1(0)$ and $b_1(0)$ are directional
derivatives in the direction of gradients. Since the boundary
$\bdS_i \cap S$ is defined equivalently both by $G_i(x)=0$ and
$g_i(x)=0$ near $u$,
the functions $G_i(x)$ and $g_i(x)$ must have parallel
gradients in the same direction at $u \in \bdS_i \cap S$. So
$a_1(0)/b_1(0)>0$ and hence $w(u)>0$. Obviously $w(x)>0$ for
interior points $x$ of $S_i$ in $S$.
\eproof

\medskip

The above lemma shows the product $g_i(x)w(x)$ has negative
definite Hessian on $S$. Unfortunately, $w(x)$ might not be a
polynomial. However, we can use polynomials to approximate $w(x)$
and its derivatives. Thus we need an improved version of {\it
Stone-Weierstrass Approximation Theorem} which shows the density
of polynomials in the space $C^k(\Omega)$ for a bounded open set $\Omega$.
Define the norm in $C^k(\Omega)$ as
\[
\| f \|_{C^k(\Omega)}:= \max_{x \in \Omega} \, \max_{\af \in \N^n,
0\leq |\af| \leq k}\, \{ |D^\af f(x) | \}.
\]

\begin{prop}
\label{thm:StonWeis}
Suppose $f \in C^k(\Omega) $ is supported in a bounded open set
$\Omega$ in $R^n$. For any $\eps>0$,  there exists a polynomial
$h$ such that $\| f- h\|_{C^k(\Omega)} <\eps$.
\end{prop}

The basic idea for proving this theorem is that
$C^{\infty}(\Omega)$ is dense in $ C^k(\Omega)$,
which contains $f$,
and then polynomials are dense in $C^{\infty}(\Omega)$.
The proof is straightforward, for example, it is an exercise
in Hirsch \cite[Chapter~2]{Hirsch}.
Thus we omit the proof here.

\begin{prop} \label{prop:approxp}
Assume $S_i=\{x\in \re^n: \, g_i(x) \geq 0\}$
is extendable poscurv-convex with respect to $S$. Then there
exists a polynomial $h_i(x)$ positive on $S$ such that the product
$p_i(x):=g_i(x)h_i(x)$ has negative definite Hessian on $S$.
\end{prop}

\noindent {\it Proof} \,
Let $T_i\supseteq S$ be a compact convex set such that
$\pt T_i$ is nonsingular
and $\pt T_i\cap S = \bdS_i \cap S$.
Then apply
Proposition~\ref{prop:ConDefFun} to $T_i$ and get a concave
defining function $G_i(x)$ for $T_i$ with negative definite
Hessian on $T_i \supseteq S$. Lemma~\ref{lem:gG} shows $w(x)=
\frac{G_i(x)}{g_i(x)} \in C^2(S)$ is  positive on $S$. So $w(x)
\in C^2(U) $ for some bounded open set $U$ containing $S$. Extend
$w(x)$ to the whole space $\re^n$ such that $w(x) =0$ for all
$x\notin U$. Let $w_\eps(x)$ be the mollified  function
\[
w_\eps(x) = \int \frac{1}{\eps^n}
\eta\left(\frac{x-y}{\eps}\right) w(y) dy
\]
where $\eta(x)$ is the {\it standard mollifier} function
\[
\eta(x) = \bca a e^{\frac{1}{\|x\|^2-1}}  & \text{ if } \|x\| <1   \\
0 & \text{ if } \|x\| \geq 1    \eca.
\]
Here the constant $a$ is chosen to make $\int_{\re^n} \eta(x) dx =1$.
The function $w_\eps(x)$ is a smooth function supported in a
bounded open set $U^\prm \supseteq U \supseteq S$.
Also $w(x)$ and
$w_\eps(x)$ are both twice differentiable on $S$, and $\|
w_\eps(x) - w(x) \|_{C^2(S)}$ can be made  arbitrarily small by
sending $\eps \to 0$.

Note that $G_i = g_i w$ is a concave function such that $-\nabla^2
G_i(x) \succ 0$ for all $x\in S$. Obviously
\[
\nabla^2 G_i(x) = w(x) \nabla^2 g_i(x) + \nabla g_i(x) \nabla
w(x)^T+ \nabla w(x)^T \nabla g_i(x) + g_i(x) \nabla^2 w(x).
\]
By Proposition~\ref{thm:StonWeis}, for arbitrary $\tau >0$, there
exists a polynomial $h_i(x)$ such that
\[\| w_\eps(x) - h_i(x) \|_{C^2(\Omega)} < \tau .\]
If $\eps$ and $\tau>0$ are small
enough, then $h_i(x)$ is positive on $S$ and the product
$p_i(x)=g_i(x) h_i(x)$ has negative definite Hessian on $S$.
\eproof

A simpler result on new defining polynomials which requires less
terminology to understand is:

\begin{cor} \label{cor:approxp}
Given a polynomial $g(x)$, if $T=\{x\in \re^n: \, g(x) \geq 0\}$
is a poscurv-convex set with nonempty interior, then there is a
polynomial $p(x)$ strictly concave on $T$ satisfying $p(x)=0, \nabla p(x)
\neq 0$ for $x \in \pt T$ and $p(x)>0$ for $x$ inside of $T$.
\end{cor}
\noindent {\it Proof} \,
Obviously, $T$ is extendable poscurv-convex with respect to
itself. By Proposition~\ref{prop:approxp}, there exists a
polynomial $h(x)$ positive on $T$ such that the product polynomial
$p(x)=g(x)h(x)$ has negative definite Hessian on $T$. If $x\in \pt
T$, then $p(x)=0$. If $x$ is in the interior of $T$, then
$p(x)>0$.
By an argument similar to that for part
(i) of Lemma~\ref{lem:gG},
$\nabla p(x)$ does not vanish on the boundary
$\pt T=\{x\in \re^n: \, g(x) = 0\}$.
\eproof

\section{Proofs}
\label{sec:proof} \setcounter{equation}{0}

As we have seen, the projections of the sets $\hS_N$ defined by
LMI~\reff{LmiN} contain the convex set $S$, for all integers
$N\geq \max_\nu d_\nu$. We need to prove that the projection actually
equals $S$ for some sufficiently large $N$. The basic idea of the
proof of this sharpness is to apply the Convex Set Separating
Theorem to produce a linear functional which is nonnegative on $S$
and negative on a given point outside $S$. We need to prove
Schm\"{u}dgen's or Putinar's representations for such linear
functionals with uniform degree bounds.
The uniform degree bounds will be proved in \S\ref{sec:complexityPosSS},
but will be used in this section.
They are Theorems~\ref{MatSmgBd} and \ref{MatPutBd}.

\subsection{Proofs of Theorem~\ref{SeqLmiPrj} and Theorem~\ref{PutLmiLift} }
\label{sec:proofConvcaveg}

Given a unit length vector $\ell \in \re^n$,
consider the optimization problem
\begin{align*}
\ell^*:=\min_{x\in \re^n} & \quad \ell^Tx \\
s.t.  & \quad g_1(x)\geq 0, \cdots, g_m(x)\geq 0.
\end{align*}
Let $u=u(\ell) \in S$ denote the minimizer, which must exist due
to the compactness of $S$. Note that $u$ must be on the boundary
$\bdS$.

Suppose $g_i(x)$ are concave on $S$ and $S$ has non-empty
interior, i.e., there exists $\xi \in S$ such that $ g_1(\xi)>0,
\cdots, g_m(\xi) >0.$ So the Slater's condition holds, which
implies that there exist nonnegative Lagrange multipliers
$\lmd_1,\cdots, \lmd_m \geq 0$ such that
\[
\ell  = \sum_{i=1}^m \lmd_i \nabla g_i(u), \quad \lmd_i g_i(u) =0,
\, \forall i = 1,\cdots, m.
\]
So $ f_\ell(x) = \ell^Tx - \ell^* - \sum_{i=1}^m \lmd_i \nabla
g_i(x) $ is a convex function on $S$ such that $f_\ell(u)=0$ and
$\nabla f_\ell(u) =0$. Hence for all $x\in S$ we have
\[
f_\ell(x) \geq f_\ell(u) + \nabla f_\ell(u)^T (x-u) = 0.
\]
So $f_\ell(x)$ is nonnegative on $S$. We hope to find a
Schm\"{u}dgen's or Putinar's representation of $f_\ell$ in terms
of polynomials $g_1(x),\cdots,g_m(x)$. But we want the
representation to have a uniform degree bound. Indeed, validating
the lifted construction in \S \ref{sec:LMIlift} amounts to proving
that there is a $N$, such that for all $\|\ell\|=1$ the
polynomials in the resulting representation of $f_\ell(x)$ have
degree at most $2N$.

\medskip

\begin{lem} \label{LinSosBd}
Use the above notations. Suppose $S$ has non empty interior and
its defining polynomials $g_i(x)$ are concave on $S$. Suppose
either $-\nabla^2g_i(x)$ is SOS or $-\nabla^2g_i(u) \succ 0$ for
all $u \in \pt S_i \cap \bdS$. Then for every unit length vector
$\ell$ we have the representation
\[
f_\ell(x) =  \sum_{i=1}^m \lmd_i \ (x-u)^T F_i(u,x) (x-u)
\]
where $u$ is the minimizer,  $\lambda_i \geq 0$ are the Lagrange
multipliers, and $F_i(u,x)$ is SOS in $x$ or such that
\be \label{eq:Fbd} \dt I_n \preceq F_i(u,x) \preceq M I_n, \quad
\forall\, x \in\ S \ee for some positive constants $M > \dt >0$
which are independent of $\ell$.
\end{lem}

\medskip

\noindent {\it Proof} \,
Since $f_\ell(u)=0$ and $\nabla f_\ell(u)=0$, we get
\[
f_\ell(x) = \sum_{i=1}^m \lmd_i \, (x-u)
\underbrace{\left(\int_0^1 \int_0^t\,
  -\nabla^2 g_i(u+s(x-u)) \, dsdt \right)}_{F_i(u,x)} (x-u).
\]
Let $J(u)=\{ 1\leq i\leq m:\, g_i(u)=0\}$ be the index set of
active constraints. For $i\notin J(u)$, the Lagrange multiplier
$\lmd_i=0$, so we can choose $F_i(u,x)$ to be the zero matrix
which is of course SOS. Note that for all $i \in J(u)$, $u\in
\bdS_i \cap \bdS$.

\medskip

If $-\nabla^2 g_i(x)$ is SOS, then $F_i(u,x)$ is also SOS in $x$
by Lemma~\ref{lem:SosInt}. If $-\nabla^2 g_i(x)$ is not SOS but
positive definite on the boundary $\bdS_i\cap \bdS$, then $F_i(u,x)$
must be positive definite on $S$. To see this point, we first show
that $F_i(u,x)$ is positive semidefinite. For any $u\in \bdS_i\cap
\bdS,\, x\in S$, the line segment $\{u+s(x-u):\, 0\leq s\leq 1\}$ is
contained in $S$ and $g_i$ are concave, so
\[
F_i(u,x) = \int_0^1 \int_0^t\,  -\nabla^2 g_i(u+s(x-u)) \, dsdt
\succeq 0.
\]
Second, we show $F_i(u,x)$ is positive definite. Suppose for some
vector $\xi$,
\[
\xi^TF_i(u,x)\xi = \int_0^1 \int_0^t \xi^T
\left(-\nabla^2g_i(u+s(x-u))\right)\xi ds dt =0.
\]
By the concavity of $g_i$, we must have
\[
\xi^T \left(-\nabla^2g_i(u+s(x-u))\right)\xi=0, \quad \forall\, s
\in [0,1].
\]
Choose $s=0$ in the above, then $-\nabla^2 g_i(u) \succ 0$ implies
$\xi =0$. Hence $F_i(u,x)\succ 0$ for all $x\in S$ and $u \in
\bdiS$.
%
%
Now we need show $F_i(u,x)$ satisfies the
inequality~\reff{eq:Fbd}. Obviously, by definition, as a function
of $u$ and $x$, $F_i(u,x)$ is continuous in $u\in \bdS_i\cap \bdS$
and $x\in S$. And $F_i(u,x)$ is positive definite for all $u\in
\bdS_i \cap S$ and $x\in S$. Since the minimum eigenvalue is a
continuous function of the matrix, the existence of constants
$M>\dt>0$
independent of $\ell$ 
is due to the compactness of $S$.
\eproof

\medskip

\begin{theorem} \label{thm:SoS}
Assume polynomials $g_i(x)$ satisfy the hypotheses of
  Theorem~\ref{SeqLmiPrj}.
Then there exists a finite integer $N$ such that for every vector
$\ell$ with $\|\ell\|=1$
\[
f_\ell(x) = \sum_{\nu\in\{0,1\}^m} \sig_\nu(x) g_1^{\nu_1}(x)
\cdots  g_m^{\nu_m}(x)
\]
where $\sig_\nu(x)$ are SOS polynomials with degree
\[
\deg(\sig_\nu g_1^{\nu_1}  \cdots  g_m^{\nu_m}) \leq 2N.
\]
Furthermore, if the archimedean condition on the $g_i$
  holds, then
$f_\ell(x)$ has the representation
\[
f_\ell(x) = \sig_0(x) + \sig_1(x) g_1(x) + \cdots  + \sig_m(x)
g_m(x)
\]
with degree bounds such that  $\deg(\sig_i g_i) \leq 2N$.
\end{theorem}

\noindent {\it Proof} \,
Let $N$ be the maximum integer such that $ \max_i \deg(g_i) \leq
2(N-1) \leq \Omega(\frac{M}{\dt})$ where function $\Omega(\cdot)$
is given by Theorem~\ref{MatSmgBd}, and $M, \delta$ are  given by
\reff{eq:Fbd}. Fix an arbitrary vector $\ell$ and let $u$ be the
minimizer of $\ell^Tx$ on $\pt S$. By Lemma~\ref{LinSosBd}
\[
 f_\ell(x) = \sum_{i=1}^m \lmd_i (x-u)^T F_i(u,x) (x-u)
 \]
holds for matrix polynomials $F_i(u,x)$ which are either SOS in
$x$ or such that
\[
\dt I_n \preceq F_i(u,x) \preceq M I_n
\]
with some positive constants $M >\dt >0$ which are independent of
$\ell$. Let $K=\{ 1\leq i\leq m:\, F_i(u,x) \text{ is SOS } \}$.
If $i\in K$, then $(x-u)^T F_i(u,x) (x-u)$ is an SOS polynomial of
degree at most $\deg(g_i)$. If $i \notin K$, by
Theorem~\ref{MatSmgBd}, there exist SOS matrix polynomials
$G_\nu^{(i)}(x)$ such that
\[
F_i(u,x) = \sum_{\nu\in\{0,1\}^m} G_\nu^{(i)}(x) g_1^{\nu_1}(x)
\cdots  g_m^{\nu_m}(x)
\]
with degree bounds $ \deg(G_\nu^{(i)} g_1^{\nu_1}  \cdots
g_m^{\nu_m}) \leq 2N-2$. Now let
\begin{align*}
\sig_0(x)  & = \sum_{i\in K} \lmd_i  (x-u)^T F_i(u,x) (x-u)
+ \sum_{i\notin K} \lmd_i  (x-u)^T G_0^{(i)}(x) (x-u) \\
\sig_\nu(x) &= \sum_{i\notin K} \lmd_i  (x-u)^T G_\nu^{(i)}(x)
(x-u), \qquad \text{ if } \nu \ne 0.
\end{align*}
which must also be SOS polynomials such that $\deg(\sig_\nu
g_1^{\nu_1}  \cdots  g_m^{\nu_m} ) \leq 2N$. So we have the
Schm\"{u}dgen representation for $f_\ell$ with uniform (in $\ell$)
degree bounds.

Similarly, Putinar's representation for $f_\ell$  with uniform
degree bounds follows from Theorem~\ref{MatPutBd}.
\eproof

\medskip

Now we are able to complete the proofs of Theorem~\ref{SeqLmiPrj} and \ref{PutLmiLift}.
The basic idea for the proof is as follows.
Theorem~\ref{thm:SoS} essentially guarantees that
the so-called S-BDNR and PP-BDNR (under archimedean condition) properties
mentioned in Section~\ref{sec:LMIlift} hold for $S$,
which implies the S-BDR and PP-BDR properties
also hold for $S$, and thus the results in \cite{Las06}
can be applied to validate the exactness of the constructed LMIs
for Theorems~\ref{SeqLmiPrj} and \ref{PutLmiLift}.

\bigskip
\noindent
{\it Proof of Theorem~\ref{SeqLmiPrj}}\,
First, we prove the S-BDNR property holds for $S$
for integer $N$ claimed by Theorem~\ref{thm:SoS}.
Let $a^Tx+b$ be nonnegative on $S$
and $u$ be a minimizer of $a^Tx+b$ on $S$.
Since $S$ has nonempty interior, the Slater's
condition holds, that is, there exist Lagrange multipliers
$\lmd_1,\cdots,\lmd_m\geq 0$ such that
\[
a = \sum_{i=1}^m \lmd_i \nabla g_i(u), \quad \lmd_i g_i(u) =
0,\quad \forall \, i =1 ,\cdots, m.
\]
Applying Theorem~\ref{thm:SoS} for $\ell = a$ and $\ell^* = a^Tu$,
we get the representation
\[
a^Tx - a^Tu - \sum_i \lmd_i g_i(x) =
\sum_{\nu\in\{0,1\}^m} \sig_\nu^\prm(x) g_1^{\nu_1}(x) \cdots
g_m^{\nu_m}(x)
\]
for some SOS polynomials $\sig_\nu^\prm$ with degree bounds
\[
\deg(\sig_\nu^\prm g_1^{\nu_1}  \cdots  g_m^{\nu_m}) \leq 2N.
\]
Or equivalently we have the identity ( note that $a^Tu + b \geq 0 $)
\be \nn
a^Tx + b =
\sum_{\nu\in\{0,1\}^m} \sig_\nu(x) g_1^{\nu_1}(x)
\cdots g_m^{\nu_m}(x)
\ee
for some SOS polynomials $\sig_\nu$ with the
same degree bounds as $\sig_\nu^\prm$.
Thus the S-BDNR property holds for $S$,
and so does S-BDR.

Second, we prove that the LMI~\reff{LmiN} $\hat S$ constructed in \S \ref{sec:LMIlift}
is a lifted LMI of $S$ for integer $N$ claimed by Theorem~\ref{thm:SoS}.
Since the S-BDNR property implies the S-BDR property for $S$,
Theorem~2 in Lasserre \cite{Las06} can be applied to
validate the exactness of the lifted LMI~\reff{LmiN}.
For the convenience of readers, we give the direct proof here,
because the proof is short and
the approach will be used in proving
Theorems~\ref{thm:mainRep} and \ref{thm:posCurv}
(these theorems can not be shown by only proving the S-BDR or PP-BDR property,
since their lifted LMIs are not purely based
on Schm\"{u}dgen's or Putinar's representation).

Obviously, the set $S$ is contained in the projection down of each
LMI $\hS_N$ defined by \reff{LmiN}.  We show they are equal. Otherwise, in
pursuit of a contradiction, suppose there exists a vector $(\hat
x, \hat y)$ in $\hS$  such that $\hat x$ is not in the convex set
$S$. By the Hahn-Banach Separation Theorem, there must exist a
vector $\ell$ of unit norm such that
\begin{align} \label{eq:minell}
\ell^T \hat x< \ell^* :=\min_{x\in S}\, \ell^Tx.
\end{align}
So $\ell^Tx-\ell^*$ is nonnegative on $S$.
By the S-BDNR property, we have the identity
\be \label{eq:felliden}
\ell^Tx - \ell^* =
\sum_{\nu\in\{0,1\}^m} \sig_\nu(x) g_1^{\nu_1}(x)
\cdots g_m^{\nu_m}(x)
\ee
for some SOS polynomials $\sig_\nu$ with degree bounds
$ \deg(\sig_\nu g_1^{\nu_1}  \cdots  g_m^{\nu_m}) \leq 2N$.
Since $\sig_\nu$ is SOS, we can write $\sig_\nu (x) = [x^{d-d_\nu}]^T W_\nu [x^{d-d_\nu}]$
for some symmetric matrix $W_\nu \succeq 0$.
Now in identity \reff{eq:felliden},
replace each monomial $x^\af$ by $\hat y_\af$, we get
\[
\ell^T \hat x - \ell^* =
\sum_{\nu\in\{0,1\}^m}  Trace\Big( W_\nu \cdot
\big(\sum_{ 0\leq |\af|\leq 2N} A_\af^\nu \hat y_\af  \big)
\Big) \geq 0,
\]
which  contradicts \reff{eq:minell}.
\eproof

\bigskip
\noindent {\it Proof of Theorem~\ref{PutLmiLift}} \,
We can prove Theorem~\ref{PutLmiLift} in a very similar manner
to the proof   above of Theorem~\ref{SeqLmiPrj}.
Here we only show the distinctive parts.

First,
by the archimedean condition and Theorem~\ref{thm:SoS},
we can
prove the PP-BDNR property holds for $S$
for integer $N$ claimed by Theorem~\ref{thm:SoS}, that is,
for any affine polynomial $a^Tx+b$ nonnegative on $S$,
we have the identity (we denote $g_0(x)=1$)
\be \nn
a^Tx + b = \sum_{k=0}^m \sig_k(x) g_k(x)
\ee
for some SOS polynomials $\sig_i$ with degree bounds
$\deg(\sig_i g_i) \leq 2N$.
This implies the PP-BDR property also holds for $S$.

Second, since PP-BDR property holds for $S$,
Theorem~2 in Lasserre \cite{Las06}
(same argument as above)
can be applied to
validate the exactness of the lifted LMI~\reff{putLMI}.
Here we directly give the proof by contradiction.
Follow the same contradiction approach
we have done in the proof of Theorem~\ref{SeqLmiPrj}.
Suppose there exists $(\hat x, \hat y)$ in $\hS$
such that $\hat x \notin S$.
Then there must exist $\ell \in \re^n$ of unit norm and $\ell^*$ such that
\be \nn
\ell^T \hat x < \ell^* :=\min_{x\in S}\, \ell^Tx.
\ee
So $\ell^Tx-\ell^*$ is nonnegative on $S$.
Then the PP-BDNR property implies the identity
\[
\ell^Tx - \ell^*  = \sum_{k=0}^m \sig_k(x) g_k(x)
\]
for some SOS polynomials $\sig_i$ with degree bounds
$\deg(\sig_i g_i) \leq 2N$.
Since $\sig_i$ is SOS, we can write $\sig_i (x) = [x^{d-d_\nu}]^T W_i [x^{d-d_\nu}]$
for some symmetric matrix $W_i \succeq 0$.
By substituting $\hat y_\af$ for each $x^\af$ in the above identity, we get
\[
\ell^T \hat x - \ell^* =
\sum_{k=0}^m  Trace\Big( W_i \cdot
\big(\sum_{ 0\leq |\af|\leq 2N} A_\af^{(k)} \hat y_\af  \big)
\Big) \geq 0,
\]
which results in the contradiction $ 0 > \ell^T \hat x - \ell^* \geq 0$.
\eproof

\subsection{Proof of Theorem~\ref{thm:conPDLH}}

In this subsection, we assume $S$ is convex, compact and has
nonempty interior, and $g_i(x)$ are concave on $S$. Then Slater's
condition holds and the Lagrange multipliers exist for the linear
objective $\ell^Tx$. For unit length vectors $\ell\in \re^n$, let
$f_\ell(x), \ell^*, u, \lmd_i$ be defined as before.

\begin{lem}
Assume the PDLH condition holds, then there exist constants
$M>\dt>0$ such that
\[
 \dt I_n \preceq
 \underbrace{\int_0^1 \int_0^t \left(
  -\sum_{i=1}^m \lmd_i \nabla^2 g_i(u+s(x-u)) \right)  dsdt}_{L(u,x)}
 \preceq M I_n
\]
for every unit length vector $\ell$.
\end{lem}
\noindent {\it Proof} \,
Let $\xi\in S$ be a fixed interior point. Note that
\[
f_\ell(x) = \ell^Tx-\ell^* - \sum_{i=1}^m \lmd_i g_i(x) \geq 0,
\quad \forall \, x\in S.
\]
Choose $x=\xi$ in the above, then we have
\[
0\leq \lmd_i \leq \frac{\ell^Tx-\ell^*}{g_i(\xi)} =
\frac{\ell^Tx-\ell^T u }{g_i(\xi)} \leq \frac{D}{g_i(\xi)}
\]
where $D$ is the diameter of $S$. So $ \max_i\lmd_i \leq
\frac{D}{\min_i g_i(\xi)}$. Thus $\lmd_i$ are uniformly bounded
and hence the existence of $M$ is obvious. Since $g_i(x)$ are
concave on $S$, we have that $L(u,x)$ must be positive
semidefinite on $S$. We need to show $\dt$ exists. Otherwise, in
pursuit of a  contradiction, suppose we have a sequence
$\{\ell^{(k)}\}$, $\{u^{(k)}\}$, $\{x^{(k)}\}$, and
$\{\lmd^{(k)}\}$ such that $\lmd_{\min}(L(u^{(k)},x^{(k)})) \to
0$. Since $\{\ell^{(k)}\}, \{u^{(k)}\}, \{x^{(k)}\},
\{\lmd^{(k)}\}$ are all bounded, without loss of generality, we
can assume
\[
\ell^{(k)} \to \hat \ell, \qquad  u^{(k)}\to \hat u, \qquad
x^{(k)}\to \hat x, \qquad \lmd^{(k)}\to \hat \lmd.
\]
The limit $\hat \ell$ also has unit length, $\hat u$ is the
minimizer of $\hat \ell^Tx$ on $S$ and $\hat \lmd_i$ are the
corresponding Lagrange multipliers. That the limit
$L(u^{(k)},x^{(k)})$ is singular implies there exists $ 0\ne
\zeta\in \re^n$ such that
\[
\int_0^1 \int_0^t \zeta^T \left(
  -\sum_{i=1}^m \hat \lmd_i \nabla^2 g_i(\hat u+s(\hat x-\hat u)) \right) \zeta  dsdt = 0.
\]
Then we must get
\[
\zeta^T \left(  -\sum_{i=1}^m  \hat\lmd_i \nabla^2 g_i(\hat
u+s(\hat x-\hat u)) \right) \zeta  = 0, \quad \forall\, 0\leq s
\leq 1.
\]
Choose $s=0$ in the above. But the PDLH condition implies $\zeta
=0$, which is a contradiction.
\eproof

\begin{theorem} \label{thm:PDLH}
Assume the hypotheses of Theorem~\ref{thm:conPDLH}. Then there
exists a finite integer $N$ such that for every unit length vector
$\ell$ with $\|\ell\|=1$
\[
f_\ell(x) = \sum_{\nu\in\{0,1\}^m} \sig_\nu(x) g_1^{\nu_1}(x)
\cdots  g_m^{\nu_m}(x)
\]
where $\sig_\nu(x)$ are SOS polynomials with degree
\[
\deg(\sig_\nu g_1^{\nu_1}  \cdots  g_m^{\nu_m}) \leq 2N.
\]
Furthermore, if the archimedean condition holds for polynomials $g_1,\cdots, g_m$, then
$f_\ell(x)$ has the representation
\[
f_\ell(x) = \sig_0(x) + \sig_1(x) g_1(x) + \cdots  + \sig_m(x)
g_m(x)
\]
with degree bounds such that  $\deg(\sig_i g_i) \leq 2N$.
\end{theorem}

\noindent {\it Proof} \,
Let $N$ be the maximum integer such that $2(N-1) \leq
\Omega(\frac{M}{\dt})$ where the function $\Omega(\cdot)$ is given by
Theorem~\ref{MatSmgBd}, and $M, \delta$ are given by the
preceding lemma. Fix an arbitrary unit length vector $\ell$. Since
$f_\ell(u)=0$ and $\nabla f_\ell(u)=0$, we get
\[
 f_\ell(x) =  (x-u)^T
 \underbrace{ \int_0^1 \int_0^t
  \left( -\sum_{i=1}^m \lmd_i \nabla^2 g_i(u+s(x-u))
  \right) dsdt}_{L(u,x)}
 (x-u).
 \]
with $ \dt I_n \preceq L(u,x)\preceq M I_n $. By
Theorem~\ref{MatSmgBd}, there exist SOS matrix polynomials
$G_\nu^{(i)}(x)$ such that
\[
L(u,x) = \sum_{\nu\in\{0,1\}^m} G_\nu(x) g_1^{\nu_1}(x) \cdots
g_m^{\nu_m}(x)
\]
with degree bounds $ \deg(G_\nu g_1^{\nu_1}  \cdots  g_m^{\nu_m})
\leq 2N-2$. Now let $ \sig_\nu(x) = (x-u)^T G_\nu(x) (x-u)$, which
must also be SOS polynomials such that $\deg(\sig_\nu g_1^{\nu_1}
\cdots  g_m^{\nu_m} ) \leq 2N$, then the first part of the theorem
holds.

\medskip
The second part of the theorem can proved by applying
Theorem~\ref{MatPutBd} in a similar way.
\eproof

\medskip

\noindent {\it Proof of Theorem~\ref{thm:conPDLH}}\,
We claim that LMI~\reff{LmiN} is an SDP representation of
$S$ when $N$ is sufficiently large.
The proof is very similar to what we have done
for proving Theorem~\ref{SeqLmiPrj}.

First, as we have done in proving Theorem~\ref{SeqLmiPrj},
Theorem~\ref{thm:PDLH} can be applied
to show the S-BDNR property, which implies
the S-BDR property holds for $S$.
So, Theorem~2 in Lasserre \cite{Las06}
validates the exactness of the lifted LMI~\reff{LmiN}.
Here we give the direct proof by contradiction,
which is very similar to what we have done for Theorem~\ref{SeqLmiPrj}.
Here we only give the distinctive parts.
Suppose there exists $(\hat x, \hat y)$ in $\hS$
such that $\hat x \notin S$.
Then there
must exist $\ell \in \re^n$ of unit norm and $\ell^*$ such that
\be \nn
\ell^T \hat x < \ell^* :=\min_{x\in S}\, \ell^Tx.
\ee
So $\ell^Tx-\ell^*$ is nonnegative on $S$.
The S-BDNR property implies that there exists an $N>0$ such that
\[
\ell^Tx - \ell^* =  \sum_{\nu\in\{0,1\}^m}
\sig_\nu(x) g_1^{\nu_1}(x) \cdots  g_m^{\nu_m}(x) .
\]
for some SOS polynomials $\sig_\nu$ with degree bounds
$\deg(\sig_\nu(x) g_1^{\nu_1}(x) \cdots  g_m^{\nu_m}(x)) \leq 2N$.
So we can write $\sig_\nu (x)
= [x^{d-d_\nu}]^T W_\nu [x^{d-d_\nu}]$
for some symmetric matrix $W_\nu \succeq 0$.
In the above identity,
similar to what we have done in the proof of Theorem~\ref{SeqLmiPrj},
replacing each monomial $x^\af$ by $\hat y_\af$ results in
the contradiction $ 0 > \ell^T \hat x - \ell^* \geq 0$.

Furthermore, if the archimedean condition holds for polynomials $g_1,\cdots, g_m$,
then we can similarly prove that
the LMI~\reff{putLMI}
is the lifted LMI for S when $N$ is big enough,
as we have done for proving Theorem~\ref{PutLmiLift}.
\eproof

\subsection{Proof of Theorem~\ref{thm:mainRep} }

In this subsection, we no longer assume the defining polynomials
$g_i(x)$ are concave on $S$ but only that they are quasi-concave.
The set $S$ is still assumed to be convex, compact and have
nonempty interior. The key point of our proof is to find a
different set of concave polynomials defining the same convex set
$S$.

Now we return to the proof of Theorem \ref{thm:mainRep}. Assume
the hypotheses of  Theorem~\ref{thm:mainRep} hold. If $-\nabla^2
g_i(x)$ is SOS in $x$, let $p_i(x)=g_i(x)$ which is obviously
concave. If $g_i$ is strictly quasi-concave on $S$,
Proposition~\ref{prop:nicefunc} implies that we can find new
defining polynomials $p_i$ that have negative definite Hessian on
$S$. So in some open set $U$ containing $S$, we have
\[
\underbrace{\{x\in\re^n:\, p_1(x) \geq 0 , \cdots, p_m(x) \geq 0
\}}_{P} \cap U = S.
\]
We should mention that the set $P$ might not coincide with $S$,
since it might be possible that for some point $v$ far away from
$S$ such that $p_i(v) \geq 0$ for all $1\leq i\leq m$. Let $y_\af
= \int x^\af d \mu(x)$ be the $\af$-th moment. Write $p_i(x) =
\sum_{\af} p^{(i)}_\af x^\af$. If $\supp(\mu) \subseteq S$, then
\[
\sum_{\af} p^{(1)}_\af y_\af \geq 0, \quad \cdots, \quad
\sum_{\af} p^{(m)}_\af y_\af \geq 0.
\]
Therefore, the set $S$ is contained in the projection of the set
$\widehat{S}_N$ of solutions to the following refined LMI
\begin{align}
\label{refLMI} \left.\baray{rl} \forall \, \nu \in \{0,1\}^m,
\quad A_0^{\nu} +
\underset{ 0 < |\af|\leq 2N}{\sum} A_\af^{\nu} y_\af & \succeq 0 \\
L_{p_1}(y) \geq 0, \cdots, L_{p_m}(y) &\geq 0 \\
x_1 = y_{e_1},\, \cdots, \, x_n=y_{e_n},\, y_0 & = 1
\earay\right\}.
\end{align}
Here $L_{p_i}(y) = \sum_{\af} p^{(i)}_\af y_\af$ and symmetric
matrices $A_\af^{\nu}$ are the same as those in LMI~\reff{LmiN}.
Note the $A_\alpha^\nu$ are determined by the $g_i$.
So \reff{refLMI} uses both $p_i$ and $g_i$.

\medskip

Our goal is to prove \reff{refLMI} is a lifted LMI for $S$ for
sufficiently large $N$, thereby validating
Theorem~\ref{thm:mainRep}. For this purpose, we need a lemma
similar to Lemma~\ref{LinSosBd}. For arbitrary unit length vector
$\ell$, consider optimization
\begin{align*}
\ell^*:=\min_{x\in U} & \quad  \ell^Tx \\
s.t. & \quad  p_1(x) \geq 0, \cdots, p_m(x)\geq 0,
\end{align*}
which is the same as to minimize $\ell^Tx$ on $S$. Let $u=u(\ell)$
be the minimizer, whose existence is guaranteed by the compactness
of $S$. Note that $S$ has an interior point $\xi\in S$, i.e.,
$g_i(\xi)>0$ for all $i=1,\cdots,m$. By our construction, $p_i(x)=
h_i(x) g_i(x)$ for some polynomials $h_i(x)$ which are positive on
$S$. So $\xi\in S$ is also an interior point for the new defining
polynomials $p_1(x),\cdots,p_m(x)$, and hence the Slater's
condition holds for the constraints $p_1(x)\geq
0,\cdots,p_m(x)\geq 0$. Therefore there exist Lagrange multipliers
$\lmd =[\, \lmd_1 \, \cdots \, \lmd_m\,]\geq 0$ such that the
function
\[
\tilde f_\ell(x) := \ell^Tx - \ell^* - \sum_{i=1}^m \lmd_i p_i(x)
\]
is a nonnegative convex function on $S$ such that $\tilde
f_\ell(u) =0,\, \nabla \tilde f_\ell(u)=0$. Note that $-\nabla^2
p_i(x)$ is either SOS or positive definite on $S$.

\medskip

\begin{lem} \label{refLinSosBd}
Let $p_i(x), \tilde f_\ell(x), \ell^*,  \lmd_i, u$ be defined as
above. Then we have the representation
\[
\tilde f_\ell(x) = \sum_{i=1}^m \lmd_i (x-u)^T F_i(u,x) (x-u)
\]
where the symmetric matrix polynomial $F_i(u,x)$ is either SOS or
such that
\[
\dt I_n \preceq F_i(u,x) \preceq M I_n, \quad \forall\, x \in S
\]
for some positive constants $M>\dt >0$ which are independent of
$\ell$.
\end{lem}

\noindent {\it Proof} \,
Since $\tilde f_\ell(u) =0,\, \nabla \tilde f_\ell(u)=0$, we have
\[
\tilde f_\ell(u) = \sum_{i=1}^m \lmd_i (x-u)^T
\underbrace{\left(-\int_0^1 \int_0^t \nabla^2
p_i(x+s(x-u)ds\,dt\right)}_{F_i(u,x)} (x-u).
\]
If $p_i(x)$ is sos-concave, $F_i(u,x)$ is SOS in $x$ by
Lemma~\ref{lem:SosInt}. If $p_i(x)$ is strictly concave on $S$,
then we can prove $F_i(u,x)$ is positive definite on $S$. Apply
the same argument in the proof for Lemma~\ref{LinSosBd}.
\eproof

\begin{theorem} \label{thm:refSoS}
Assume polynomials $g_i(x)$ satisfy the hypotheses of
Theorem~\ref{thm:mainRep}. Then there exists a finite integer $N$
such that for every vector $\ell$ with $\|\ell\|=1$
\[
\tilde f_\ell(x) = \sum_{\nu\in\{0,1\}^m} \sig_\nu(x)
g_1^{\nu_1}(x) \cdots  g_m^{\nu_m}(x)
\]
where $\sig_\nu(x)$ are sums of squares of polynomials with degree
\[
\deg(\sig_\nu g_1^{\nu_1}  \cdots  g_m^{\nu_m}) \leq 2N.
\]
Furthermore, if the archimedean condition on the $g_i$
  holds, then
$\tilde f_\ell(x)$ has the representation
\[
\tilde f_\ell(x) = \sig_0(x) + \sig_1(x) g_1(x) + \cdots  +
\sig_m(x) g_m(x)
\]
with degree bounds such that  $\deg(\sig_i g_i (x)) \leq 2N$.
\end{theorem}

\medskip
\noindent {\it Proof} \,
The proof is almost the same as for Theorem~\ref{thm:SoS}.
Just follow the argument for proving Theorem~\ref{thm:SoS}.
The only differences are replacing
$f_\ell(x)$ by $\tilde f_\ell(x)$
and then applying Lemma~\ref{refLinSosBd}
instead of Lemma~\ref{LinSosBd}.
\eproof

\medskip

\noindent{\it Proof of Theorem~\ref{thm:mainRep}}\,
If $-\nabla^2 g_i(x)$ is SOS in $x$, let $p_i(x)=g_i(x)$ which is
obviously concave. If $g_i$ is strictly quasi-concave on $S$, let
$p_i(x)$ be the new defining polynomials for $S_i$ given by
Proposition~\ref{prop:nicefunc}, which have negative definite
Hessian on $S$. For some small open set $U$ containing $S$, the
convex set $S$ is equivalently defined as
\[
S = \{x\in U:\, p_1(x) \geq 0, \cdots, p_m(x) \geq 0 \}.
\]
As we have seen earlier, $S$ is contained in the projection of
LMI~\reff{refLMI}. We claim that this projection is sharp for $N$
given by Theorem~\ref{thm:refSoS}. The proof is very similar to
the proof for Theorem~\ref{SeqLmiPrj}.

\smallskip

Otherwise, seeking a contradiction, suppose there exists a vector
$(\hat x, \hat y)$ in $\hat S_N$ such that $\hat x \notin S$. By
the Hahn-Banach Separation Theorem, there must exist a vector
$\ell$ of unit length such that
\be  \label{eq:ldown}
\ell^T \hat x< \ell^* :=\min_{x\in S}\, \ell^Tx.
\ee
Let $u\in S$ be the minimizer of $\ell^T x$ on $S$, which must be
on the boundary $\pt S$. Note that $ p_1(x), \cdots, p_m(x)$ are
concave polynomials, and $S$ has nonempty interior. Since $p_i(x)
= h_i(x) g_i(x)$ for $h_i(x)$ positive on $S$, the new equivalent
constraints $ p_1(x)\geq 0, \cdots, p_m(x)\geq 0$ also have
nonempty interior. Thus the Slater's condition holds and hence
there exist Lagrange multipliers $\lmd_1,\cdots,\lmd_m\geq 0$ such
that
\[
\ell =  \sum_{i=1}^m  \lmd_i \nabla p_i(u), \quad \lmd_i p_i(u) =0
,\ \forall\, i=1,\cdots,m.
\]
By Lemma~\ref{refLinSosBd} and Theorem~\ref{thm:refSoS}, we get
\[
\tilde f_\ell(x)=\ell^Tx - \ell^* - \sum_{i=1}^m \lmd_i p_i(x) =
\sum_{\nu\in\{0,1\}^m} \sig_\nu (x) g_1^{\nu_1}(x) \cdots
g_m^{\nu_m}(x)
\]
for some SOS polynomials $\sig_\nu$ with
$\deg(\sig_\nu g_1^{\nu_1}  \cdots  g_m^{\nu_m}) \leq 2N. $
So we have the identity
\[
\ell^Tx - \ell^*  = \sum_{i=1}^m \lmd_i p_i(x) +
\sum_{\nu\in\{0,1\}^m} \sig_\nu(x) g_1^{\nu_1}(x) \cdots
g_m^{\nu_m}(x).
\]
We can write $\sig_\nu (x)
= [x^{d-d_\nu}]^T W_\nu [x^{d-d_\nu}]$
for some symmetric matrix $W_\nu \succeq 0$.
In the above identity,
replacing each monomial $x^\af$ by $\hat y_\af$,
we get
\[
\ell^T \hat x - \ell^*  = \sum_{i=1}^m \lmd_i L_{p_i}(\hat y) +
\sum_{\nu\in\{0,1\}^m} Trace\Big(
W_\nu \cdot \big( \sum_{0 \leq |\af|\leq 2N} A_\af^{\nu} \hat y_\af \big)
\Big) \geq 0,
\]
which contradicts \reff{eq:ldown}.
\eproof

\medskip

\noindent
{\it Remark:}
In LMI (~\ref{refLMI}), we use all the
products $g^\nu(x) = g_1^{\nu_1}(x) \cdots g_m^{\nu_m}(x)$ for all
$\nu\in \{0,1\}^m$ which results an exponential size of LMI. As we
did in the end of Section~\ref{sec:LMIlift}, the set $S$ is also
the projection of the following set
\begin{align} \label{refPutLMI}
\left.\baray{rl} \forall \, 0\leq k \leq m,  \quad A_0^{(k)} +
\underset{ 0 < |\af|\leq 2N}{\sum} A_\af^{(k)} y_\af & \succeq 0 \\
L_{p_1}(y) \geq 0, \cdots, L_{p_m}(y) &\geq 0 \\
x_1 = y_{e_1},\, \cdots, \, x_n=y_{e_n},\, y_0 & = 1
\earay\right\}
\end{align}
where symmetric matrices $A_\af^{(k)}$ are defined in
LMI~\reff{putLMI}. If the archimedean condition holds, we can
similarly prove \reff{refPutLMI} is a lifted LMI for $S$ when $N$
is sufficiently large, as we did in the above proof.

\subsection{Proof of Theorems~\ref{thm:setSimple} and \ref{thm:posCurv} }

In the remarks after Theorem~\ref{thm:posCurv}, we mentioned that
Theorems~\ref{thm:setSimple} can be implied by Theorem~\ref{thm:posCurv}.
So we only need to prove Theorem~\ref{thm:posCurv}.

\begin{lem} \label{lem:US}
Let $S$ be as in Theorem~\ref{thm:posCurv}.
Then there exists an open set containing $S$
and polynomials $p_1(x),\cdots,p_m(x)$
which either are sos-concave or have negative definite Hessian on $U$
such that
\[
S = \{x\in U:\, p_1(x) \geq 0, \cdots, p_m(x) \geq 0 \}.
\]
\end{lem}
\noindent {\it Proof} \,
If $S_i = \{x\in\re^n: g_i(x)\}$ is sos-convex, choose $p_i(x) = g_i(x)$.
If $S_i$ is extendable poscurv-convex with respect to $S$,
by Proposition~\ref{prop:approxp},
there exists a polynomial $h_i(x)$ positive on $S$
such that the product $g_i(x)h_i(x)$ has negative definite Hessian on $S$,
then choose $p_i(x) = g_i(x) h_i(x)$.
Since $S$ is compact, we can choose an open set $U \supset S$
small enough to make the lemma true.
\eproof
\bigskip
\noindent {\it Proof of Theorem~\ref{thm:posCurv}}\,
The proof is almost the same as the one for Theorem~\ref{thm:mainRep}.
We follow the approach there, and only list the distinctive parts here.
Let $U$ and polynomials $p_i(x)$ be given by Lemma~\ref{lem:US}.
Then define LMI~\reff{refLMI} using both $p_i$ and
$g_i$, and $S$ is contained in the projection of this LMI.
Then we claim \reff{refLMI} is a lifted LMI for $S$
for $N$ given by Theorem~\ref{thm:refSoS}.

Similarly, we prove this by contradiction.
Suppose there exists a vector
$(\hat x, \hat y)$ in $\hat S_N$ such that $\hat x \notin S$. By
the Hahn-Banach Separation Theorem, there must exist a vector
$\ell$ of unit length and $\ell^*$ such that
\be \nn
\ell^T \hat x- \ell^* < 0, \qquad
\ell^Tx - \ell^* \geq 0 , \, \forall \, x\in S.
\ee
Let $u$ be a minimizer of $\ell^Tx$ over $S$.
Then Slater's condition implies the first order optimality condition holds at $u$
for the set of defining polynomials $\{x\in U:  p_1(x) \geq 0, \cdots, p_m(x) \geq 0 \}$.
Let $\lmd_i \geq 0$ be the corresponding Lagrange multipliers and
\[
\tilde f_\ell (x) := \ell^Tx - \ell^* - \sum_{i=1}^m \lmd_i p_i(x).
\]
Similarly, by Theorem~\ref{thm:refSoS}, we can get the identity
\[
\ell^Tx - \ell^*  = \sum_{i=1}^m \lmd_i p_i(x) +
\sum_{\nu\in\{0,1\}^m} \sig_\nu(x) g_1^{\nu_1}(x) \cdots
g_m^{\nu_m}(x)
\]
for some SOS polynomials $\sig_\nu$ with degree bounds
$\deg(\sig_\nu g_1^{\nu_1}  \cdots  g_m^{\nu_m}) \leq 2N$.
We can also write $\sig_\nu (x)
= [x^{d-d_\nu}]^T W_\nu [x^{d-d_\nu}]$
for some symmetric matrix $W_\nu \succeq 0$.
Similar to what we have done in the proof for Theorem~\ref{thm:mainRep},
a contradiction to $\ell^T\hat x - \ell^* \geq 0$ can be obtained
by replacing each monomial $x^\af$ by $\hat y_\af$ in the above identity.
\eproof

\bigskip
\section{Appendix: The complexity of the matrix Positivstellensatz}
\label{sec:complexityPosSS} \setcounter{equation}{0}

Throughout this section, we only need assume $S = \{x\in \re^n:\,
g_1(x)\geq 0, \cdots, g_m(x) \geq 0\}$ is compact. We do not need
either $g_i(x)$ is concave or $S$ is convex. Without loss of
generality, assume $S \subset (-1,1)^n$, otherwise do some
coordinate transformation.

\medskip

Suppose we have a symmetric matrix polynomial $F(x) \in
\re[x]^{r\times r}$ which is positive definite on $S$. Our goal is
to give  a Positivstellensatz representation of $F(x)$ in terms of
defining polynomials $g_1(x), \cdots, g_m(x)$ with bounds on the
degrees of the representing polynomials. When $r=1$, that is,
$F(x)$ are scalar polynomials, Schm\"{u}dgen's Positivstellensatz
\cite{Smg} says that $F(x)$ has the representation
\[
F(x) =  \sum_{\nu \in \{0,1\}^m } \sig_\nu(x) g_1^{\nu_1}(x)
\cdots g_m^{\nu_m}(x)
\]
for some SOS polynomials $\sig_\nu$. Furthermore, if the
archimedean condition holds, Putinar's Positivstellensatz
\cite{Putinar} says that $F(x)$ has the representation
\[
F(x) = \sig_0(x) + \sig_1(x) g_1(x)+ \cdots + \sig_m(x)  g_m(x)
\]
for some SOS polynomials $\sig_i$.

\medskip

These representation results can be generalized to the case $r>1$.
Schm\"{u}dgen's matrix Positivstellensatz says that there
exist symmetric SOS matrix polynomials $G_\nu(x) \in
\re[x]^{r\times r}$ such that
\[
F(x) =  \sum_{\nu \in \{0,1\}^m } g_1^{\nu_1}(x) \cdots
g_m^{\nu_m}(x) G_\nu(x).
\]
Under the archimedean condition, Putinar's matrix
Positivstellensatz says that there exist symmetric SOS matrix
polynomials $G_i(x) \in \re[x]^{r\times r}$ such that
\[
F(x) = G_0(x) +  g_1(x) G_1(x)+ \cdots + g_m(x) G_m(x).
\]
We refer to \cite{SH} for these representations of positive
definite matrix polynomials. The goal of this section is to give
degree bounds for $G_\nu(x)$ in these representations.

\subsection{Schm\"{u}dgen's matrix Positivstellensatz}

For a scalar polynomial $f(x) = \sum_{\af}  f_\af x^\af$, its norm
$\|f\|$ is defined to be
\be \label{def:normf}
\| f \| = \max_\af \left\{ |f_\af|  \frac{\af_1 ! \cdots \af_n ! }
{ |\af| !} \right\}.
\ee
For a matrix polynomial $F(x) =  \sum_{\af}  F_\af x^\af$, its norm is defined to be
\be  \label{def:normF}
\| F \| = \max_\af \left\{ \|F_\af\|_2 \frac{\af_1 ! \cdots \af_n ! }
{ |\af| !} \right\} =
\max_\af \left\{ \frac{ \|\pt_{x_1}^{\af_1} \cdots \pt_{x_n}^{\af_n } F(0) \|_2 }
{ |\af| !} \right\}  .
\ee
Here $\|A\|_2$ denotes the maximum singular value of matrix $A$.

\begin{lem} \label{lem:polya}
Suppose polynomials $g_i(x)$ are scaled such that for some
$\veps>0$
\[
S = \{x\in (-1+\veps, 1-\veps)^n:\ \ g_1(x) \geq 0, \cdots, g_m(x)
  \geq 0,\ \
\sum_{i=1}^m g_i(x) < 2n\veps \}.
\]
Define new polynomials
\begin{align*}
p_1 &= 1-\veps+x_1, \cdots, p_n = 1-\veps+x_n \\
p_{n+1} &= 1-\veps-x_1, \cdots, p_{2n} = 1-\veps-x_n \\
p_{2n+1} &= g_1, \cdots, p_{2n+m} = g_m,\ \; \  p_{2n+m+1}
=2n\veps-(g_1+\cdots+g_m).
\end{align*}
Then there exists an integer $c>0$  depending only on the
polynomials $g_1,\cdots,g_m$ such that
for every symmetric matrix polynomial $F(x) \succ 0$ on $ S$
can be written
\begin{align} \label{eqf}
F = \sum_{|\af| \leq N }  p_1^{\af_1} \cdots
p_{2n+m+1}^{\af_{2n+m+1}} F_\af
\end{align}
with constant symmetric matrices $F_\af \succ 0$ and $N \leq
\Theta\left(\frac{\|F\|}{F^*}\right)$.
Here
\[
d =\deg(F(x)), \quad
F^* : = \min_{x\in S} \lmd_{\min} (F(x)), \quad
\Theta(s): = c d^2
\left( 1+\left( d^2n^{d} s\right)^c\right).
\]
\end{lem}

\noindent {\it Proof} \,
The proof is very similar to the one of Lemma~9 by Schweighofer \cite{Schw}
and uses a famous theorem of Poly\'a.
The only difference is the scalar polynomials $f(x)$ in \cite{Schw}
are replaced by matrix polynomials $F(x)$.
We follow the approach for proving Lemma~9 in \cite{Schw}.
Without loss of generality, assume $\|F\|=1$.
Introduce new variables
$y = (y_1,\cdots, y_{2n+m+1})$.
Define the homomorphism
\[
\varphi: \re[y] \to \re[x] :\,\,
y_i \mapsto p_i.
\]
Then $\varphi(y_1+\cdots+y_{2n+m+1}) = 2n$
and hence $y_1+\cdots+y_{2n+m+1} - 2n \in \ker (\varphi)$.
By Hilbert's basis theorem, there exist polynomials
$r_1, \cdots, r_t$ so that
\[
\ker (\varphi) = \langle  y_1+\cdots+y_{2n+m+1} - 2n, r_1, \cdots, r_t \rangle.
\]
Then define new sets
\begin{align*}
\Dt  &:= \{ y \in \re_+^{2n+m+1}:\, y_1+\cdots+y_{2n+m+1} = 2n\}, \\
Z& :=  \{ y \in \Dt:\, r_1(y) = \cdots = r_t(y) =0 \}.
\end{align*}
The properties listed below hold,
which are essentially Claim~1 and Claim~2 of \cite{Schw}.
\bdes
\item [(P1)] The linear map
\[
\ell: \re^{2n+m+1}   \to \re^n:  \,\,
(y_1,\cdots,y_{2n+m+1}) \mapsto \half (y_1-y_{n+1}, \cdots, y_n-y_{2n})
\]
induces a bijection $\ell\big|_Z : Z \to S$.

\item [(P2)] There exists an integer $d_0\geq 1$ and a $d_0$-form $R_0\in \ker\varphi$
such that $R_0\geq 0$ on $\Dt$ and $Z=\{y\in \Dt:\,\, R_0(y) =0 \}$.

\edes
By Lojasiewicz's inequality (Corollary~2.6.7 in \cite{BCR}), there exist
integers $c_0,c_1\geq 1$ such that
\[
\mathbf{(E1)}\hspace{2cm}
\mbox{dist}(y,Z)^{c_0} \leq c_1 R_0(y), \quad \forall \,\, y\in \Dt.
\hspace{4cm}
\]
Define new constants
\[
c_2 := 2^{c_0+1}c_1 \sqrt{2n}, \, \qquad
c_3:=c_2 (2n)^{d_0} \|R_0\|, \, \qquad
c_4 := (2n)^{d_0}.
\]
Then choose $c>0$ big enough so that
\[
d_0^2 (1+c_4a +c_3a^{c_0+1}) \leq c(1+a^c),\quad
\forall\,\, a \in [0, \infty).
\]

Now we write $F(x) = F_0(x)+\cdots + F_d(x)$ with
$F_k$ being matrix $k$-forms
(homogeneous matrix polynomials of degree $k$).
Set $d_1:= \max\{d,d_0\}$ and
\[
P(y):=\sum_{k=0}^d P_k(y)
\left(\frac{y_1+\cdots+y_{2n+m+1}}{2n}\right)^{d_1-k},
\qquad
P_k(y) :=  F_k\left(\half (y_1-n_{n+1}), \cdots,
\half (y_n - y_{2n}) \right).
\]
Then $P$ is a $d_1$-form such that
\[
\varphi(P) = F, \qquad \qquad P(y) = F (\ell(y)) \,\,\forall y\in \Dt.
\]
So we can see (P1) implies
\[
\mathbf{(E2)} \hspace{1cm}
\min \left\{ \lmd_{\min} \Big(P(y)\Big) :\, y\in Z \right\} =
\min \left\{ \lmd_{\min} (F(x)) :\, x\in S \right\} = F^* >0.
\hspace{2cm}
\]
Next, define the $d_1$-form $R$ as follows
\[
R(y):= R_0(y) \cdot
\left(\frac{y_1+\cdots+y_{2n+m+1}}{2n}\right)^{d_1-d_0}.
\]
By equations (17) and (18) of \cite{Schw}, we know
$
\|R\| \leq \frac{1}{(2n)^{d_1-d_0}} \| R_0\|
$
and
\[
\mathbf{(E3)} \hspace{3cm} R(y) = R_0(y),\quad\forall\, y\in \Dt.
\hspace{4cm}
\]

\bigskip

Now we claim the property listed below holds
\bdes
\item [(P3)]  For all $y,y^\prm \in \Dt$, it holds
\[
| \lmd_{\min}\Big(P(y)\Big) - \lmd_{\min}\Big(P(y^\prm)\Big) |
\leq \| P(y) - P(y^\prm) \|_2 \leq \sqrt{n}d^2 n^{d-1} \| y - y^\prm\|.
\]
\edes
The first inequality of (P3) can be obtained by noting the fact that
for any two symmetric matrices $A,B$ it holds
\begin{align*}
\lmd_{\min}(A) &= \min_{ \|\xi\|=1}  \xi^T A \xi =
 \min_{ \|\xi\|=1} \left( \xi^T B \xi  + \xi^T(A-B)\xi \right)
\leq \lmd_{\min}(B)  + \|A-B\|_2, \\
\lmd_{\min}(B) & = \min_{ \|\xi\|=1}  \xi^T B \xi =
\min_{ \|\xi\|=1} \left(  \xi^T A \xi  + \xi^T(B-A)\xi \right)
\leq \lmd_{\min}(A)  + \|A-B\|_2.
\end{align*}
The second inequality of (P3) is a consequence of
Claim~3 in \cite{Schw}, because
\[
\| P(y) - P(y^\prm) \|_2 = \sup_{ \|\xi\|=1}
| \xi^T P(y)\xi - \xi^T P(y^\prm) \xi |
\]
and Claim~3 in \cite{Schw} can be applied to
the scalar polynomials $\xi^T P(y)\xi $.

\bigskip

For those $y, y^\prm \in \Dt$ satisfying
$\lmd_{\min}\Big(P(y)\Big) \leq \half F^*$ and
$\lmd_{\min}\Big(P(y^\prm)\Big) \geq F^*$, by (P3), we have
\[
\| y - y^\prm\|  \geq \frac{F^*}{2\sqrt{n}d^2 n^{d-1}}
\geq \frac{F^*}{2d^2 n^d}.
\]
Therefore, Property~(E2) implies that,
for all $y\in\Dt$ with $\lmd_{\min}\Big(P(y)\Big)
\leq \half F^*$, it holds
\[
\mbox{dist}(y,Z) \geq \frac{F^*}{2d^2n^d},
\]
and hence Properties (E2) and (E3) imply,
for all $y\in\Dt$ with $\lmd_{\min}\Big(P(y)\Big) \leq \half F^*$,
it holds
\[
\left(\frac{F^*}{2d^2n^d}\right)^{c_0} \leq c_1 R(y).
\]
In (P3), if we choose $y^\prm$ to be a minimizer of
$\lmd_{\min}\Big(P(y)\Big)$
on $Z$, then for all $y\in \Dt$ we have
\[
|\lmd_{\min}\Big(P(y)\Big) - F^* | \leq \mbox{diameter}
(\Dt) \sqrt{n}d^2n^{d-1}
\leq 2\sqrt{2n}d^2n^d
\]
which obviously implies
\[
\mathbf{(E4)} \hspace{3cm}
\lmd_{\min}\Big(P(y)\Big) \geq F^* -   2\sqrt{2n}d^2n^d, \qquad
\forall \, y\in \Dt.
\hspace{4cm}
\]
Let $\lmd := c_2 d^2 n^d \left( \frac{d^2 n^d}{F^*}\right)^{c_0}$
and define a new set
\begin{align*}
\Dt_1 := \{ y\in \re_+^{2n+m+1}:\, y_1+\cdots+y_{2n+m+1} = 1\}.
\end{align*}
Now we claim that
\[
\mathbf{(E5)} \hspace{4cm}
\lmd_{\min}  \Big( P(y) \Big)  + \lmd R(y)    \geq
\half F^*,  \forall \, y \in \Dt
\hspace{4cm}
\]
where $r$ is the dimension of the matrix polynomial $F(x)$.
Now we prove (E5).
Obviously, (E5) holds for those $y\in \Dt$ with
$\lmd_{\min}\left(  P(y) \right)   \geq \half F^*$.
We only need to verify (E5) for those
$y \in \Dt$ with $\lmd_{\min}\left(P(y) \right)  \leq \half F^*$.
The choice of $\lmd$ shows
\[
\lmd R(y) \geq \frac{c_2}{c_12^{c_0}} d^2n^d
\]
and hence (E4) implies
\[
\lmd_{\min}\Big(P(y)\Big)   + \lmd R(y)  \geq
F^* -   2\sqrt{2n}d^2n^d + \frac{c_2}{c_12^{c_0}} d^2n^d
= F^* \geq \half F^*.
\]
Therefore we obtain that (by concavity of function $\lmd_{\min}\big(\cdot\big)$)
\begin{equation*}
\lmd_{\min}\Big(P(y)  + \lmd R(y) I_r\Big)  \geq
\lmd_{\min}\Big(P(y)\Big) + \lmd R(y) \geq
\half F^*, \quad \forall \, y \in \Dt
\end{equation*}
which by homogeneity implies
\[
\lmd_{\min}\Big(P(y)  + \lmd R(y) I_r\Big)
\geq  \frac{F^*}{2(2n)^{d_1}}, \quad \forall \, y \in \Dt_1.
\]
Then Theorem~3 from Scherer and Hol \cite{SH} guarantees that the product
\[
Q(y):= \Big(P(y)+\lmd R(y) \cdot I_r\Big)
\cdot \left( \frac{y_1+\cdots+y_{2n+m+1} }{2n} \right)^N
\]
has positive definite matrix coefficients for all
\[
N > \frac{d_1(d_1-1)\|P+\lmd R \cdot I_r \|}
{2 \frac{F^*}{2(2n)^{d_1}}} -d_1
= d_1(d_1-1) (2n)^{d_1} \frac{\|P+\lmd R \cdot I_r \|} {F^*} - d_1.
\]
If $N$ is chosen to be the smallest integer in the above, then
\[
\deg(Q) \leq cd^2 \left( 1+  \left(\frac{d^2n^d}{F^*}\right)^c\right),
\]
as is shown at the end of the proof of Lemma~9 in \cite{Schw}.
Since $F(x)= \varphi (P(x)) = \varphi (Q(x))$,
we have proved $F(x)$ can be represented like \reff{eqf}
with $F_\af \succ 0$.
\eproof

Now we arrive at our theorem giving degree bounds.
\begin{theorem}
\label{MatSmgBd} If matrix polynomial $F(x) \succeq \dt I \succ 0$
for all $x$ in  a compact subset $S$ of $ \re^n$, then
\[
F(x) = \sum_{ \nu \in \{0,1\}^m  } g_1^{\nu_1} \cdots g_m^{\nu_m}
G_\nu(x)
\]
where $G_\nu(x)$ are SOS matrix polynomials such that
\[
\deg(g_1^{\nu_1} \cdots g_m^{\nu_m} G_\nu) \leq
\Omega\left(\frac{\|F\|}{\dt}\right) :=\kappa \cdot
\Theta\left(\frac{\|F\|}{\dt}\right)
\]
where $\kappa$ is a constant depending only on the polynomials
$g_i(x)$.
\end{theorem}

\noindent {\it Proof} \,
Again take $S \subset (-1+\veps, 1-\veps)^n $. Then the
polynomials $p_1,\cdots, p_{2n}, p_{2n+m+1}$ in the preceding
lemma are positive on $S$. By Schm\"{u}dgen's Positivstellensatz,
for every $i \in \{1,\cdots, 2n, 2n+m+1\}$ we have \be
\label{eq:pRep} p_i(x)  = \sum_{ \nu \in \{0,1\}^m }
\sig_\nu^{(i)} g_1^{\nu_1} \cdots g_m^{\nu_m} \ee with
$\sig_\nu^{(i)}(x)$ being SOS polynomials.  Let
\[
\kappa = \max \left\{ \deg(\sig_\nu^{(i)} g_1^{\nu_1} \cdots
g_m^{\nu_m}):\, i=1,\cdots,2n, 2n+m+1,\, \nu \in \{0,1\}^m \right
\}.
\]
Now into  identity \reff{eqf} we plug the representation
\reff{eq:pRep} for $p_1(x), \cdots, p_{2n}(x), p_{2n+m+1}(x)$,
then we obtain the conclusion in the theorem.
\eproof

\subsection{Putinar's matrix Positivstellensatz}

\begin{lem}\label{cubelift}
Suppose $g_i(x)$ are scaled such that $g_i(x)\le 1$ on $[-1,1]^n$.
Then there exist constants $c_0,c_1,c_2>0$ with the property:
\bcen
\begin{minipage}{0.9\textwidth}
For all symmetric matrix polynomials $F(x)\in\re[x]^{r\times r}$
of degree $d$ such that $F(x) \succeq \dt I_r$ for all $x\in S$,
if we set
\begin{equation*}
L:=d^2n^{d-1}\frac{\|F\|}{\dt},\qquad
\lmd:=c_1d^2n^{d-1}\|F\|L^{c_2}
\end{equation*}
and if $k\in\N$ satisfies
\begin{equation*}
2k+1\ge c_0(1+L^{c_0}),
\end{equation*}
then the inequality
\begin{equation*}
F(x)-\lmd\sum_{i=1}^m(g_i(x)-1)^{2k}g_i(x)  I_r \succeq
\frac{\dt}2 I_r
\end{equation*}
holds on $[-1,1]^n$ and hence on the unit ball
$B(0,1)=\{x\in\re^n:\, \|x\| \leq 1\}$.
\end{minipage}
\ecen
\end{lem}

\noindent {\it Proof} \,
Apply Lemma~13 in \cite{NS} to the polynomial $ f^{(\xi)}(x): =
\xi^TF(x)\xi$ where $\xi \in \re^r$ is a unit length vector. Note
that the minimum value of $f^{(\xi)}(x)$ is at least $\dt$ and
$\|f^{(\xi)}(x)\|$ is at most $\|F\|$. For $L,\lmd, k$ given in the
lemma, we have for all $x \in [-1,1]^n$
\[
\xi^TF(x)\xi-\lmd\sum_{i=1}^m(g_i(x)-1)^{2k}g_i(x) \geq
\frac{\dt}2,\quad \forall\, \xi \in \re^r,\, \|\xi\|=1
\]
which implies
\[
F(x)-\lmd\sum_{i=1}^m(g_i(x)-1)^{2k}g_i(x)  I_r \succeq
\frac{\dt}2 I_r
\]
for all $x \in [-1,1]^n$.
\eproof

\begin{theorem}
\label{MatPutBd} Assume the archimedean condition holds for the
$g_i$.
If $F(x)$ is a symmetric matrix polynomial of degree $d$ such that
$F(x) \succeq \dt I \succ 0$ for all $x\in S \subset
  \re^n$,
then
\begin{align}\label{matputSOS}
F(x) = G_0(x) +  g_1(x) G_1(x) + \cdots + g_m(x) G_m(x)
\end{align}
where $G_i(x)$ are SOS matrix polynomials such that
\[
\deg(G_0), \   \deg(g_1 G_1),\cdots, \deg(g_m G_m)\  \leq \ c
\left(d^2 n^d \frac{\|F\|}{\dt}\right)^c
\]
for some constant $c$ depending only on the polynomials $g_i(x)$.
\end{theorem}
\noindent {\it Proof} \,
The proof is almost the same as for Theorem~6 in \cite{NS}. By the
archimedean condition, we can assume \be \label{eq:ac} R -
\sum_{i=1}^n x_i^2 = s_0(x) + s_1(x)g_1(x) + \cdots + s_m(x)g_m(x)
\ee for SOS polynomials $s_i(x)$. Without loss of generality, we
can assume $R=1$, because otherwise we can apply some coordinate
transformation. Let
\[
d_1 = \max_i (\deg(s_ig_i)), \qquad d_2 = 1 + \max_i (\deg(g_i) ).
\]
First, apply Lemma~\ref{cubelift} to find constants $L,\, \lmd, k$
such that
\[
\tilde F(x) :=F(x)-\lmd\sum_{i=1}^m(g_i(x)-1)^{2k}g_i(x) I_r
\succeq  \frac{\dt}2 I_r,\quad x\in [-1,1]^n.
\]
By (33) in the proof of Theorem~6 in \cite{NS}, we have for every
$\|\xi\|=1$
\[
\| \xi^T \tilde F(x)\xi  \|  \leq   \| \xi^T F(x)\xi  \| + \lmd
d_2^{2k+1}
\]
which implies
\[
\| \tilde F \| \leq  \| F \|  + \lmd d_2^{2k+1}.
\]
By (34) in the proof of Theorem~6 in \cite{NS}, we get
\[
\deg ( \tilde F(x) ) \leq \max\{d, (2k+1)d_2, 1\} := d_h.
\]
Then we apply Theorem~\ref{MatSmgBd} to $\tilde F(x)$ on the unit
ball $B(0,1)$. So there exists some constant $c_3>0$ such that \be
\label{eq:MatSosBall} \tilde F(x) = H_0(x) + (1 - \sum_{i=1}^n
x_i^2) H_1(x) \ee for some SOS matrix polynomials $H_i(x)$ with
degree
\[
\deg(H_0) \leq k_h, \quad 2+\deg(H_1) \leq k_h \qquad where \ \
k_h:= c_3 d_h^2 \left( 1 + d_h^2 n^{d_h} \frac{2 \| \tilde
F\|}{\dt} \right)^{c_3}.
\]
By (36) in the proof of Theorem~6 in \cite{NS}, we know for some
constant $c_4>0$
\[
k_h \leq c_4 \exp(c_4 d_h) \left( \frac{2 \| \tilde F\|}{\dt}
\right)^{c_3}.
\]
By (37) in the proof of Theorem~6 in \cite{NS}, we have for some
constant $c_5>0$
\[
 \frac{2 \| \tilde F\|}{\dt} \leq \exp( L^{c_5}).
\]
Combine \reff{eq:ac} and \reff{eq:MatSosBall} to get
\[
\tilde F(x) = \underbrace{ H_0(x) + s_0(x) H_1(x)}_{G_0(x)} +
g_1(x) \underbrace{s_1(x) H_1(x)}_{G_1(x)} + \cdots + g_m(x)
\underbrace{s_m(x) H_1(x)}_{G_m(x)}.
\]
Now we can estimate $\deg(g_i G_i)$ by following the proof for
Theorem~6 in \cite{NS}. The techniques are exactly same. Finally
we can obtain the degree bound in \reff{matputSOS} for some
constant $c>0$.
\eproof

\section{Conclusions}  \label{sec:cncl}

This paper studies the SDP representation of convex sets.
Obviously, for a set $S$ to be SDP representable,
necessary conditions are that $S$ must be convex
and semialgebraic.
It is not known if these conditions are
also sufficient, but the main
contribution of this paper is to give some additional conditions
which are sufficient.
Given $S=\{x\in\re^n:\,g_1(x)\geq 0, \cdots, g_m(x)\geq 0\}$
which is  convex, compact and
has nonempty interior, we have proved $S$ is SDP representable
in either of the following cases:
(i) All $g_i(x)$
are  concave on $S$, and the PDLH condition holds;
(ii) Each $S_i$ is either sos-convex or extendable
poscurv-convex with respect to $S$.
%

The key to our proofs is to find and prove a well-behaved
Schm\"{u}dgen or Putinar's representation
for the affine polynomial $\ell^Tx-\ell^*$ nonnegative on $S$,
that is, to find conditions for the S-BDR property
or the PP-BDR property to hold for affine polynomials.
When $\ell^Tx-\ell^*$ is nonnegative on
$S$ and equals  zero at $u$ in $S$,
we can not directly apply
Schm\"{u}dgen or Putinar's Positivstellensatz to prove the
existence of the representation. However, we should mention that
it is possible to prove the existence of Schm\"{u}dgen or
Putinar's representations for $\ell^Tx-\ell^*$ by
applying the  representation of nonnegative polynomials in
Marshall \cite{Mar06} and Scheiderer \cite{Sched03,Sched05}. But
the degrees of these representations depend on the choice of
$\ell$ and we can not get a uniform degree bound from these
papers. So we were motivated to use Hessians of defining
polynomials to get the degree bound independent of $\ell$.
One interesting future work is to get
the SDP representability of $S$ by using methods in Marshall
\cite{Mar06} and Scheiderer \cite{Sched03,Sched05}.

The main result of this paper is that
if the boundary of every $S_i$ is either sos-convex or
extendable poscurv-convex with respect to $S$,
then the compact convex set $S$ is semidefinite representable.
We point out that the condition of extendable poscurv-convexity
does not require much more than that
the boundary $\bdS_i \cap \bdS$ is nonsingular and has positive curvature.
In the follow-up paper \cite{HN2} to this one,
the authors have proved a stronger result:
if for every $i$ either $S_i$ is sos-convex or $\bdS_i \cap \bdS$ is
positively curved and nonsingular,
then the compact convex set $S$ is SDP representable;
this is based on constructions of the lifted LMIs and
theorems we have proved in this paper.
Since a necessary condition for a set $S$ to be convex is that
its boundary $\bdS$ can have only nonnegative curvature
(under some nonsingularity assumption on $\bdS$),
we can see that the sufficient conditions
of semidefinite representability
given in this paper are not far away from the necessary conditions.

\bigskip
\noindent {\bf Acknowledgements.}  We are grateful to the IMA in
Minneapolis for hosting us during our collaboration on this
project. J.~William Helton was partly supported
by the NSF DMS 0700758, DMS 0400794
and the Ford Motor Co.
We thank J.Lasserre and P. Parrilo for
conversations on the lifted LMIs and to J. Lasserre for sending us
an early version of \cite{Las06}. We thank S.~Baouendi,
P.~Ebenfeld, B.~Sturmfels, and R.~Williams for discussions
on polynomial approximation,
and two referees on their comments on improving this paper.
We also thank M.~Schweighofer for pointing out
a gap in the original proof of Lemma~\ref{lem:polya}.

\end{document}